\documentclass[12pt]{article}
\usepackage[utf8]{inputenc}
\usepackage[T1]{fontenc}
\usepackage[]{geometry}
\usepackage{amssymb}
\usepackage{url}
\usepackage{systeme}
\usepackage{url}
\usepackage{array}
\usepackage{authblk}
\usepackage{lscape}
\usepackage{float} 

\usepackage{hyperref}
\hypersetup{colorlinks=true,urlcolor=Blue,linkcolor=black,citecolor=Rhodamine}
	\usepackage[dvipsnames,natural,table,xcdraw]{xcolor}

\makeatletter
\newcommand*{\rom}[1]{\expandafter\@slowromancap\romannumeral #1@}
\makeatother

\newcommand{\poch}[2]{\left({#1}\right)_{#2}}
\newcommand{\pp}[1]{\left({#1}\right)}
\newcommand{\bb}[1]{\left[{#1}\right]}

\newmuskip\pFqskip
\mathchardef\pFcomma=\mathcode`, 

\newcommand*\pFq[5]{%
  \begingroup
  \begingroup\lccode`~=`,
    \lowercase{\endgroup\def~}{\pFcomma\mkern\pFqskip}%
  \mathcode`,=\string"8000
  {}_{#1}F_{#2}\biggl(\genfrac..{0pt}{}{#3}{#4}\Big|#5\biggr)%
  \endgroup
}

\usepackage{tikz}
\usepackage{amsmath}
\usepackage{tikz-cd}
\usetikzlibrary{matrix, calc, arrows}

\usepackage{stackrel}
\newcommand{\leftrarrows}{\mathrel{\raise.75ex\hbox{\oalign{%
  $\scriptstyle\leftarrow$\cr
  \vrule width0pt height.5ex$\hfil\scriptstyle\relbar$\cr}}}}
\newcommand{\lrightarrows}{\mathrel{\raise.75ex\hbox{\oalign{%
  $\scriptstyle\relbar$\hfil\cr
  $\scriptstyle\vrule width0pt height.5ex\smash\rightarrow$\cr}}}}
\newcommand{\Rrelbar}{\mathrel{\raise.75ex\hbox{\oalign{%
  $\scriptstyle\relbar$\cr
  \vrule width0pt height.5ex$\scriptstyle\relbar$}}}}
\newcommand{\longleftrightarrows}{\leftrarrows\joinrel\Rrelbar\joinrel\lrightarrows}

\makeatletter
\def\leftrightarrowsfill@{\arrowfill@\leftrarrows\Rrelbar\lrightarrows}
\newcommand{\xleftrightarrows}[2][]{\ext@arrow 3399\leftrightarrowsfill@{#1}{#2}}
\makeatother
\date{\today}

\numberwithin{equation}{section}

\usepackage{adjustbox}
\usepackage{amssymb}
\usepackage{bbold}
\usepackage{incgraph}

\usepackage{cleveref}

\crefformat{section}{\S#2#1#3} 
\crefformat{subsection}{\S#2#1#3}
\crefformat{subsubsection}{\S#2#1#3}

\begin{document}
\title{\bf Continuous $-1$ Hypergeometric\\
Orthogonal Polynomials}
\author[1]{Jonathan Pelletier}
\author[1,2]{Luc Vinet}
\author[3]{Alexei Zhedanov}
\affil[1]{Centre de recherches mathématiques, Université de Montréal, P.O. Box 6128, Centre-ville
	Station, Montréal (Québec), H3C 3J7, Canada}
\affil[2]{IVADO, 6666 Rue Saint-Urbain, Montréal (Québec), H2S 3H1, Canada}
\affil[3]{Department of Mathematics, School of Information, Renmin University of China, Beijing 100872, China}
\date{\today}
\maketitle
\hrule
\begin{abstract}
The study of $-1$ orthogonal polynomials viewed as $q\to -1$ limits of the $q$-orthogonal polynomials is pursued. This paper present the continuous polynomials part of the $-1$ analog of the $q$-Askey scheme. A compendium of the properties of all the continuous $-1$ hypergeometric polynomials and their connections is provided. 
\end{abstract}

\hrule
\section{Introduction}
In the last decade, attention has been paid to classical orthogonal polynomials that arise as limits when $q$ goes to $-1$ of polynomials of the $q$-Askey scheme \cite{L-1J,BI,CBI,CHI}. These are referred to as $-1$ polynomials. By classical, we mean that in addition to being orthogonal and hence obeying a three-term recurrence relation, they satisfy an eigenvalue equation of the form
\begin{gather}
	L \mathbf{P}_n\pp{x} = \lambda_n \mathbf{P}_n\pp{x},
\end{gather}
where $L$ is a differential or difference operator. A striking feature of $-1$ polynomials is that they are eigenfunctions of Dunkl type operators involving reflections \cite{dunkl}. A well known example of such polynomials was introduced by Bannai and Ito (Bannai-Ito polynomials) in \cite{BIbook}, but it was only recently that their eigenvalue equation was identified \cite{BI} and that their bispectrality was hence established. This eigenvalue equation was found by looking at the most general symmetric first order Dunkl shift operator, that preserves the degree of the polynomials. The Bannai-Ito and the complementary Bannai-Ito polynomials are the most general ones. Being both $q\to -1$ limits of the $q$--Racah polynomials, they sit at the top of the $-1$ scheme that is emerging. These families of orthogonal polynomials bring two important concepts in the characterization of the $-1$ polynomials. The first is the notion of spectral transformations and kernel polynomials. Indeed, these two families are kernel partners, as one can be obtained via a spectral transformation (Christoffel/Geronimus) of the other \cite{BI}. This type of connection will recur in the description of the continuous part of the $-1$ scheme. The second concept is the Leonard duality. A monic discrete family with this property is solution to a three-term difference equation in addition to the three-term recurrence relation
\begin{gather}
    x\mathbf{P}_n(x) = \mathbf{P}_{n+1}(x) + b_n \mathbf{P}_n(x) + u_n \mathbf{P}_{n-1}(x)\label{rec0},
\end{gather}
where $b_n,u_n\in \mathbb{R}$ and $u_n>0$. In the discrete $q$--Askey scheme, all polynomials enjoy this property. It is not the case for $-1$ polynomials in general. The complementary Bannai-Ito polynomials and the dual $-1$ Hahn polynomials \cite{CBI,D-1H} were found to satisfy a five-term difference equation. This brings new issues in the definition of dual families of orthogonal polynomials under the exchange of the variable and the degree. This dual operation will lead to a five-term recurrence relation, which lies beyond the scope of ordinary orthogonal polynomials (see \cite{CHIbook}). These two concepts will also appear when considering the continuous part of the $-1$ scheme.

The main purpose of this paper is to give a classification of the continuous $-1$ orthogonal polynomials and to organize them in a scheme corresponding to the $q\to -1$ limit of the continuous q-Askey tableau. Four categories of $-1$ continuous families will be introduced, with many relations connecting them. The characterization of the members of each family, including the explicit expression in terms of hypergeometric series, the three-term recurrence relation, the orthogonality relation, the Dunkl difference or differential equation and the relations to other families of polynomials in the scheme will be provided. The approach is based on the different limits and specializations of the recurrence coefficients of known families of $q$ and $-1$ orthogonal polynomials. From Favard theorem, it is given that the resulting polynomials will be orthogonal if $b_n,u_n\in \mathbb{R}$ and $u_n>0$ in (\ref{rec0}). Focusing on continuous polynomials, the corresponding part of the $-1$ scheme will be constructed. Starting with the most general polynomials with four parameters at the top of this scheme, all other families will be found by cascading down the scheme via specializations or limit processes. Going from a level to the one below, a parameter is lost, until no parameter remains.

The ($q$-) Askey scheme, with its simplicity, is very useful and practical. Recently, Koornwinder has raised in \cite{chart_AS_1} three issues with this scheme that call for possible improvement. These have to do with the completeness of the scheme, the significance (or insignificance) of certain families in the classification and the uniformity of the transformations between families. Following ideas of Verde-Star \cite{Verde_Star}, he presented the $q$--Verde--Star scheme, and more recently, the $q$--Zhedanov scheme \cite{chart_AS_2} focusing on associated algebras. The present work will limit itself to the framework of the $q$-Askey scheme to determine its $-1$ analog for continuous polynomials. The methods used to derive the $q$--Verde--Star and $q$--Zhedanov schemes are nicely systematic, their application to the $q=-1$ case will however be deferred to future work.

The paper will unfold as follows. In each section, one of the four categories of continuous $-1$ orthogonal polynomials will be presented. The most general family within this class will be characterized, and its links to the other members of the category will be established through specializations and limit processes. In Section 2, the category of $-1$ polynomials obtained from direct specialization of the big $-1$ Jacobi polynomials and the Chihara polynomials is presented. Generalized Gegenbauer and little $-1$ Jacobi polynomials are also included, as well as the Gegenbauer and special little $-1$ Jacobi polynomials. The special role played by spectral transformations in structuring the scheme is emphasized in this section. In the next section, we present the $-1$ polynomials that cascade from the continuous Bannai-Ito polynomials. They include two families: the continuous Bannai-Ito and the continuous $-1$ Hahn polynomials. The first one is seen to be connected to the big $-1$ Jacobi and the second, presumably new, is characterized and connected to a $q\to -1$ limit of the continuous $q$ Hahn polynomials. In Section 4, the $-1$ Meixner-Pollaczek polynomials are introduced as a generalization of the generalized Hermite polynomials via a $q\to -1$ limit of the $q$-Meixner-Pollaczek polynomials. They are also connected to the continuous $-1$ Hahn polynomials of Section 3 and to the Chihara polynomials of Section 2. In Section 4, a continuous equivalent of the complementary Bannai-Ito is investigated. These polynomials do not form an orthogonal set, but some specializations of it yield orthogonal families. This leads to the inclusion of the generalized symmetric Bannai-Ito polynomials in the scheme of continuous $-1$ orthogonal polynomials. Finally, the complete scheme and a compendium of the main properties and connections of the continuous $-1$ orthogonal polynomials is provided in the appendix.

\section{Orthogonal Polynomials descending from\\ $q\to -1$ limit of the big $q$-Jacobi polynomials}
In this section, the bulk of the continuous part of the $-1$ scheme is presented. It is the component that has been the most documented to date (see \cite{L-1J}, \cite{B-1J} and \cite{CHI}). The Chihara polynomials and the big $-1$ Jacobi polynomials are at the top of this category, with three parameters each. These two sets of polynomials are related by spectral transformations (Christoffel and Geronimus) with spectral parameter equal to 1 \cite{CHI}. They both correspond to $q\to -1$ limit of the same family of polynomials, the big $q$-Jacobi polynomials. From \cite{koek}, we know that the monic big and little $q$-Jacobi polynomials are related as follows
\begin{gather}
\mathbf{P}_n^{(\text{big})}\pp{x;a,b,0|q} = (aq)^n\mathbf{P}_n^{\pp{\text{little}}}\pp{\frac{x}{aq};b,a|q}
\end{gather}
and we also know that there exists a $q\to -1$ limit of the little $q$-Jacobi polynomials: the little $-1$ Jacobi polynomials. This raises a few questions. Does the little $q$-Jacobi polynomials have a second $q\to -1$ limit and, if so, are those two limits connected via the same spectral transformation? It turns out that it is the case. The recurrence relation of this dilated little q-Jacobi is as follows.
\begin{gather}
    x\mathbf{P}_n(x) = \mathbf{P}_{n+1}(x) + \pp{1-A_n+C_n} \mathbf{P}_n(x) + A_{n-1}C_n \mathbf{P}_{n-1}(x)\label{I_B-1JRR},
\end{gather}
with $\mathbf{P}_{-1}(x)=0$ and $\mathbf{P}_0(x)=1$ and where
\begin{gather}
    \begin{array}{l}
A_{n}= \frac{\left(1-b q^{n+1}\right)\left(1-ab q^{n+1}\right)}{\left(1-a b q^{2 n+1}\right)\left(1-a b q^{2 n+2}\right)}, \qquad
C_{n}= \frac{ab^2 q^{2n+1}\left(1-q^{n}\right)\left(1-a q^{n}\right)}{\left(1-a b q^{2 n}\right)\left(1-a b q^{2 n+1}\right)}.
\end{array}
\end{gather}
It is straightforward to show that using $q=-e^\varepsilon$, $a=-e^{\varepsilon\alpha}$ and $b=-e^{\varepsilon\beta}$ in (\ref{I_B-1JRR}), one obtains the same definition of the little -1 Jacobi obtained in \cite{L-1J} after taking the limit $\varepsilon\to 0$. Moreover, using $q=-e^\varepsilon$, $a=-e^{\varepsilon\pp{2\alpha+1}}$ and $b=e^{2\varepsilon\beta}$, one finds the recurrence relation (\ref{I_RGG}): 
\begin{gather}
    x\mathbf{G}_n(x) = \mathbf{G}_{n+1}(x) +  \sigma_n \mathbf{G}_{n-1}(x)\label{I_RGG},
\end{gather}
with $\mathbf{G}_{-1}(x)=0$ and $\mathbf{G}_0(x)=1$ and where
\begin{gather}
    \sigma_{2 n}=\frac{n(n+\beta)}{(2 n+\alpha+\beta)(2 n+\alpha+\beta+1)}, \quad \sigma_{2 n+1}=\frac{(n+\alpha+1)(n+\alpha+\beta+1)}{(2 n+\alpha+\beta+1)(2 n+\alpha+\beta+2)}.
\end{gather}
This set of polynomials is clearly positive definite for all $n$ if $\alpha>-1$ and $\beta>0$. The recurrence relation is identified as the one for the monic generalized Gegenbauer polynomials \cite{GG}. When the three-term recurrence relation is presented in the form (\ref{I_B-1JRR}), it suffices to use the map $A_n\to C_{n+1}$ and $C_n\to A_n$ to take the Christoffel transformation with parameter 1 \cite{CHIbook}. The generalized Gegenbauer polynomials can also be recovered using this procedure, and so the little $-1$ Jacobi polynomials are the kernel partner of the generalized Gegenbauer polynomials. The generalized Gegenbauer polynomials are also a specialization of the Chihara polynomials, and thus the structure of the three-parameter row repeats itself at the two-parameter level. This replication of the structure is due to the compatibility of the specializations and of the limits from one row to the lower one. In a diagram, this translates into the commutativity of the sub-diagram containing the two rows. This commutativity puts restrictions on the possible limits that can be used to go down one row in the scheme.
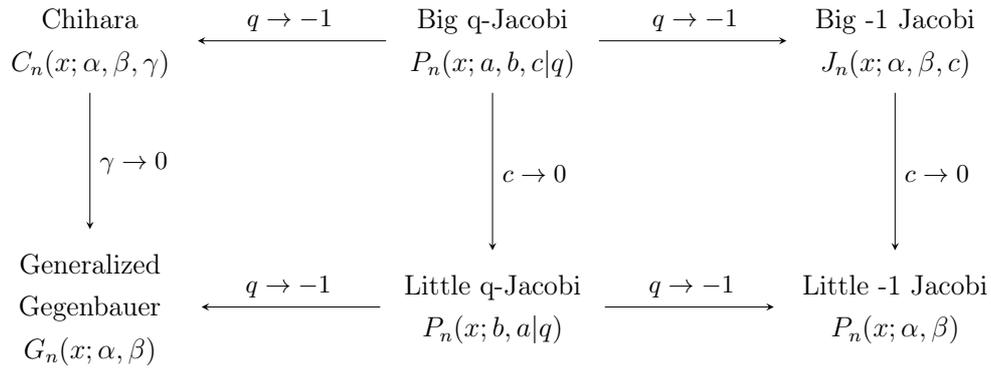
\begin{figure}[H]
\begin{center}
\adjustbox{max width=14cm}{ \begin{tikzpicture}[>=stealth,->,shorten >=2pt,looseness=.5,auto]
\matrix [matrix of math nodes,
column sep={6cm,between origins},
row sep={4cm,between origins}]
{
|(21)| \begin{tabular}{c} \( \text{Chihara} \) \\ \( C_n(x;\alpha,\beta,\gamma) \) \end{tabular} & |(22)|  \begin{tabular}{c} \( \text{Big q-Jacobi} \) \\ \( P_n(x;a,b,c|q) \) \end{tabular} & |(23)| \begin{tabular}{c} \( \text{Big -1 Jacobi} \) \\ \( J_n(x;\alpha,\beta,c) \) \end{tabular} \\
|(31)| \begin{tabular}{c} \( \begin{array}{c}
       \text{Generalized}\\\text{Gegenbauer}
      \end{array} \) \\ \( G_n(x;\alpha,\beta) \) \end{tabular} & |(32)| \begin{tabular}{c} \( \text{Little q-Jacobi} \) \\ \( P_n(x;b,a|q) \) \end{tabular} & |(33)| \begin{tabular}{c} \( \text{Little -1 Jacobi} \) \\ \( P_n(x;\alpha,\beta) \) \end{tabular} \\
};
\tikzstyle{every node}=[font=\small\itshape]
\draw (22) -- node[midway,right] {$c\to 0$}(32);
\draw (21) -- node[midway,right] {$\gamma\to 0$}(31);
\draw (23) --  node[midway,right] {$c\to 0$}(33);
\draw (22) -- node[midway,above] {$q\to-1$} (21);
\draw (32) -- node[midway,above] {$q\to-1$} (31);
\draw (22) -- node[midway,above] {$q\to-1$} (23);
\draw (32) -- node[midway,above] {$q\to-1$} (33);
\end{tikzpicture}
}
\end{center}
\caption{Correspondence between the two different $q\to -1$ limits of the big $q$-Jacobi and those of the little $q$-Jacobi}
\end{figure}
This leads to the inclusion of a third row composed of the Gegenbauer, the special little $q$-Jacobi and the special little $-1$ Jacobi polynomials. They correspond to specializations of the two-parameter row using $\alpha=-\frac{1}{2}$, $a=-1$ and $\alpha=0$ respectively. The special little $-1$ Jacobi polynomials, even if they were not specifically identified before, were introduced in \cite{L-1J} and were seen to be a solution to a Sturm-Liouville equation in addition to the Dunkl-differential equation. This leads to a factorization of the Hamiltonian operator as the square of a first order Dunkl-differential operator. The special little $-1$ Jacobi and the Gegenbauer polynomials are also kernel partner with parameter 1.

\section{-1 Orthogonal Polynomials descending from\\ the continuous Bannai-Ito polynomials}

In this section, the continuous Bannai-Ito polynomials are shown to be the parent of both the big -1 Jacobi polynomials and the continuous -1 Hahn polynomials. This last family is introduced below, and its main properties are derived through its connection to the continuous Bannai-Ito polynomials. 

The continuous Bannai-Ito polynomials were introduced in \cite{NON_SYM} as a continuous version of the Bannai-Ito polynomials. These two families are closely related (using complex parameters and orthogonality on an imaginary axis) and they also correspond to $q\to-1$ limits of the Askey-Wilson polynomials. The explicit expression of the continuous Bannai-Ito $\mathbf{Q}_n(x)$ is:\\
\begin{minipage}{\textwidth}
\begin{flalign}
    &\frac{\mathbf{Q}_{2n}\left(x;\alpha,\beta,\gamma,\delta\right)}{\pp{-2i}^{2n}} = &\label{CoBIBHR1}\end{flalign}\vspace{-0.75cm}\begin{align}\nonumber
    &\xi_{2n}\kappa_{n-1}^{(2)}\left(\frac{ix}{2}-\mathfrak{b}-\frac{1}{2}\right)\pFq{4}{3}{-n+1,n+g+2,\frac{ix}{2}+\mathfrak{b}+1,\frac{-ix}{2}+\mathfrak{b}+\frac{3}{2}}{\frac{5}{2}+\mathfrak{a}+\mathfrak{b},2+\mathfrak{b}+\mathfrak{c},2+\mathfrak{b}+\mathfrak{d}}{1} + \\
    &\kappa_n^{(1)}\pFq{4}{3}{-n,n+g+1,\frac{ix}{2}+\mathfrak{b},\frac{-ix}{2}+\mathfrak{b}+\frac{1}{2}}{\frac{3}{2}+\mathfrak{a}+\mathfrak{b},1+\mathfrak{b}+\mathfrak{c},1+\mathfrak{b}+\mathfrak{d}}{1}, \nonumber\end{align}
    \end{minipage}
    \begin{minipage}{\textwidth}
    \begin{flalign} &\frac{\mathbf{Q}_{2n+1}\left(x;\alpha,\beta,\gamma,\delta\right)}{\pp{-2i}^{2n+1}} =&\label{CoBIBHR2}\end{flalign}\vspace{-0.75cm}\begin{align}\nonumber &\kappa_{n}^{(2)}\left(\frac{ix}{2}-\mathfrak{b}-\frac{1}{2}\right)\pFq{4}{3}{-n,n+g+2,\frac{ix}{2}+\mathfrak{b}+1,\frac{-ix}{2}+\mathfrak{b}+\frac{3}{2}}{\frac{5}{2}+\mathfrak{a}+\mathfrak{b},2+\mathfrak{b}+\mathfrak{c},2+\mathfrak{b}+\mathfrak{d}}{1}+\\
    &\eta_{2n+1}\kappa_n^{(1)}\pFq{4}{3}{-n,n+g+1,\frac{ix}{2}+\mathfrak{b},\frac{-ix}{2}+\mathfrak{b}+\frac{1}{2}}{\frac{3}{2}+\mathfrak{a}+\mathfrak{b},1+\mathfrak{b}+\mathfrak{c},1+\mathfrak{b}+\mathfrak{d}}{1} \nonumber,
\end{align}
\end{minipage}
where
\begin{gather}
    g=\mathfrak{a}+\mathfrak{b}+\mathfrak{c}+\mathfrak{d}+1,\qquad \mathfrak{a}=\overline{\mathfrak{d}}=\alpha+i\beta, \qquad \mathfrak{b}=\overline{\mathfrak{c}}=\gamma+i\delta,
\end{gather}
\begin{gather}
\begin{array}{ll}
\xi_{2n}=\frac{n\pp{n+\mathfrak{c}+\mathfrak{d}+\frac{1}{2}}}{\left(2n+g\right)}, & \kappa_{n}^{(1)}=\frac{\poch{\frac{3}{2}+\mathfrak{a}+\mathfrak{b}}{n} \poch{1+\mathfrak{b}+\mathfrak{c}}{n} \poch{1+\mathfrak{b}+\mathfrak{d}}{n}}{\poch{n+g+1}{n}}, \\[0.15cm]
\eta_{2n+1}=\frac{\pp{n+\mathfrak{b}+\mathfrak{c}+1}\pp{n+\mathfrak{b}+\mathfrak{d}+1}}{\pp{2n+g+1}}, &  \kappa_{n}^{(2)}=\frac{\poch{\frac{5}{2}+\mathfrak{a}+\mathfrak{b}}{n} \poch{2+\mathfrak{b}+\mathfrak{c}}{n} \poch{2+\mathfrak{b}+\mathfrak{d}}{n}}{\poch{n+g+2}{n}}.
\end{array}
\end{gather}
Their recurrence relation reads:
\begin{gather}
    x \mathbf{Q}_{n}(x)=\mathbf{Q}_{n+1}(x)+b_{n} \mathbf{Q}_{n}(x)+u_{n} \mathbf{Q}_{n-1}(x)\label{ICBIRR},
\end{gather}
where
\begin{gather}
    b_{n}=\left\{\begin{array}{ll}
2 \beta-\frac{(n+4 \alpha+2)(\beta-\delta)}{(n+2 \alpha+2 \gamma+2)}-\frac{n(\beta+\delta)}{(n+2 \alpha+2 \gamma+1)} & n \text { even } \\
2 \beta-\frac{(n+4 \alpha+4 \gamma+3)(\beta+\delta)}{n+2 \alpha+2 \gamma+2}-\frac{(n+4 \gamma+1)(\beta-\delta)}{n+2 \alpha+2 \gamma+1} & n \text { odd }
\end{array}\right.,
\end{gather}
and 
\begin{gather}
    u_{n}=\left\{\begin{array}{ll}
\frac{n(n+4 \alpha+4 \gamma+2)\|n+2[\alpha+\gamma+i(\beta+\delta)]+1\|^{2}}{4(n+2 \alpha+2 \gamma+1)^{2}} & n \text { even } \\
\frac{(n+4 \alpha+1)(n+4 \gamma+1)\|n+2[\alpha+\gamma+i(\beta-\delta)]+1\|^{2}}{4(n+2 \alpha+2 \gamma+1)^{2}} & n \text { odd }
\end{array}\right.,
\end{gather}
with $\alpha,\beta,\gamma,\delta\in\mathbb{R}^+$.

\subsection{A Limit to the Big -1 Jacobi Polynomials}
As a limit of the Askey-Wilson polynomials with four parameters, the continuous Bannai-Ito polynomials are poised to be at the top of the -1 continuous scheme. We now establish their connection to the big $-1$ Jacobi polynomials. Let $\alpha,\beta,\gamma$ and $\delta$ be parametrized as follows
\begin{gather}
 \alpha=a_1,\quad \beta=\frac{b_1}{h},\quad \gamma=a_2,\quad \delta=\frac{b_2}{h},
\end{gather}
and denote by
\begin{gather}
    \mathbf{Q}_{n}^{(h)}(x) =\pp{\frac{h}{2b_1}}^n\mathbf{Q}_{n}\pp{\frac{2b_1 x}{h}},
\end{gather}
the renormalized polynomials resulting from the scaling $x\to\frac{2b_1 x}{h}$. Upon taking the limit $h\to0$, one obtains the polynomials $\mathbf{Q}_{n}^{(0)}(x)$ satisfying the recurrence relation
\begin{gather}
    x\mathbf{Q}_{n}^{(0)}(x)) = \mathbf{Q}_{n+1}^{(0)}(x) + \left(1-A_n-C_n\right)\mathbf{Q}_{n}^{(0)}(x) + A_{n-1}C_n \mathbf{Q}_{n-1}^{(0)}(x),
\end{gather}
\begin{minipage}{.45\textwidth}
        \begin{eqnarray}
    A_{n}=\left\{\begin{array}{ll}
\frac{\pp{1-\frac{b_2}{b_1}}\pp{n+4a_1+2}}{2n+4a_1+4a_2+4}&n\text{ even}\\[0.15cm]
\frac{\pp{1+\frac{b_2}{b_1}}\pp{n+4a_1+4a_2+3}}{2n+4a_1+4a_2+4}& n\text{ odd}
\end{array}\right.\nonumber,
\end{eqnarray}
    \end{minipage}%
    \begin{minipage}{0.55\textwidth}
       \begin{eqnarray}
   C_{n}=\left\{\begin{array}{ll}
\frac{\pp{1+\frac{b_2}{b_1}}n}{2n+4a_1+4a_2+2}&n\text{ even}\\[0.15cm]
\frac{\pp{1-\frac{b_2}{b_1}}\pp{n+4a_2+1}}{2n+4a_1+4a_2+2}& n\text{ odd}
\end{array}\right.\label{B-1JRRC}.
\end{eqnarray}
    \end{minipage}\\[0.3cm]
Those polynomials $\mathbf{Q}_{n}^{(0)}(x)$ are identified as the monic big -1 Jacobi polynomials. The continuous Bannai-Ito polynomials are thus the parent of the big -1 Jacobi polynomials, and it is appropriate to have them at the top of the -1 continuous scheme. The exact connection is
\begin{gather}
    \lim\limits_{h\to0}\pp{\frac{h}{2b_1}}^n\mathbf{Q}_{n}\pp{\frac{2b_1 x}{h};a_1,\frac{b_1}{h},a_2,\frac{b_2}{h}}=\mathbf{J}_n\pp{x;4a_1+1,4a_2+1,-\frac{b_2}{b_1}}.
\end{gather}

\subsection{Continuous -1 Hahn Polynomials}\label{C-1HP}
The monic continuous q-Hahn polynomials (with $a=c$, $b=d$ and $x\to\pp{1+q}x$ in the definition of \cite{koek}) are defined through the recurrence relation
\begin{gather}
x \mathbf{P}_{n}(x)=\mathbf{P}_{n+1}(x)+\frac{1}{2}\left[\frac{a \mathrm{e}^{i \phi}+a^{-1} \mathrm{e}^{-i \phi}}{1+q}-\left(A_{n}+C_{n}\right)\right] \mathbf{P}_{n}(x)
+\frac{1}{4} A_{n-1} C_{n} \mathbf{P}_{n-1}(x),
\end{gather}
where
\begin{gather}
    A_{n}=\frac{\left(1-a b \mathrm{e}^{2 i \phi} q^{n}\right)\left(1-a^2 q^{n}\right)\left(1-ab q^{n}\right)\left(1-a^2 b^2 q^{n-1}\right)}{a \mathrm{e}^{i \phi}\pp{1+q}\left(1-a^2 b^2 q^{2 n-1}\right)\left(1-a^2 b^2 q^{2 n}\right)},
\end{gather}
and 
\begin{gather}
    C_{n}=\frac{a \mathrm{e}^{i \phi}\left(1-q^{n}\right)\left(1-a b q^{n-1}\right)\left(1-b^2 q^{n-1}\right)\left(1-a b \mathrm{e}^{-2 i \phi} q^{n-1}\right)}{\pp{1+q}\left(1-a^2 b^2 q^{2 n-2}\right)\left(1-a^2 b^2 q^{2 n-1}\right)},
\end{gather}
with $ \mathbf{P}_{-1}(x)=0$ and $ \mathbf{P}_{0}(x)=1$. Besides the $q\to1$ limit leading to the continuous Hahn polynomials, a non-trivial $q\to -1$ limit exist and is obtained by using the parametrization
\begin{gather}
    q=-e^\varepsilon,\quad a=e^{\varepsilon\pp{2\alpha+1}},\quad b=e^{\varepsilon\pp{2\gamma+1}},\quad \phi=\frac{\pi}{2}+2\varepsilon\beta.
\end{gather}
Let $\mathbf{K}_{n}^{\pp{1}}(x)$ be the polynomials obtained when $\varepsilon\to 0$. A direct calculation of the recurrence coefficients gives the relation
\begin{gather}
    x \mathbf{K}_{n}^{\pp{1}}(x)=\mathbf{K}_{n+1}^{\pp{1}}(x)+b_{n} \mathbf{K}_{n}^{\pp{1}}(x)+u_{n} \mathbf{K}_{n-1}^{\pp{1}}(x), \label{recC-1HIn}
\end{gather}
where
\begin{gather}
    b_{n}=\left\{\begin{array}{ll}
2 \beta-2\beta\frac{n}{(n+2 \alpha+2 \gamma+1)} & n \text { even } \\
2 \beta-2\beta\frac{(n+4 \alpha+4 \gamma+3)}{n+2 \alpha+2 \gamma+2} & n \text { odd }
\end{array}\right.,
\end{gather}
and 
\begin{gather}
    u_{n}=\left\{\begin{array}{ll}
\frac{n(n+4 \alpha+4 \gamma+2)\|n+2[\alpha+\gamma+2i\beta]+1\|^{2}}{4(n+2 \alpha+2 \gamma+1)^{2}} & n \text { even } \\
\frac{(n+4 \alpha+1)(n+4 \gamma+1)\|n+2[\alpha+\gamma]+1\|^{2}}{4(n+2 \alpha+2 \gamma+1)^{2}} & n \text { odd }
\end{array}\right.,
\end{gather}
with $ \mathbf{K}_{-1}^{\pp{1}}(x)=0$ and $ \mathbf{K}_{1}^{\pp{1}}(x)=1$. It is easy to identify those continuous -1 Hahn polynomials as a direct specialization of the continuous Bannai-Ito presented in (\ref{CoBIBHR1})-(\ref{CoBIBHR2}) with $\delta=\beta$. The complete characterization (expression, difference equation and orthogonality relation) of these continuous $-1$ Hahn polynomials can be deduced from the properties of the continuous Bannai-Ito ones. In particular, their bispectrality follows from the fact that they satisfy the complex difference eigenvalue equation
\begin{gather}
    L \mathbf{K}_{n}^{(1)}(x)=\lambda_{n} \mathbf{K}_{n}^{(1)}(x), \quad \lambda_{n}=(-1)^{n}(n+2\alpha+2\gamma+3 / 2),\\[0.3cm]
L =\mathcal{A}\left(S^{+} R-I\right)
+\overline{\mathcal{A}}\left(S^{-} R-I\right)+(2\alpha+2\gamma+3 / 2) I,\\
\mathcal{A} = \frac{\pp{2 \alpha+1+i\bb{\beta- x}}\pp{2 \gamma+1+i\bb{\beta- x}}}{1-2 i x},
\end{gather}
where $\overline{\mathcal{A}}$ is the complex conjugate of $\mathcal{A}$, $S^{\pm}f(x)=f(x\pm i)$ and $Rf(x)=f(-x)$. There is in addition another limit resulting from the parametrization
\begin{gather}
    q=-e^\varepsilon,\quad a=e^{\varepsilon\pp{2\alpha+1}},\quad b=-e^{\varepsilon\pp{2\gamma+1}},\quad \phi=\frac{\pi}{2}+2\varepsilon\beta,
\end{gather}
which gives polynomials $\mathbf{K}_{n}^{\pp{2}}(x)$. These correspond to continuous Bannai-Ito polynomials with $\delta=-\beta$. This second type of continuous $-1$ Hahn polynomials is very similar to the first type. The differences between the two types do not affect the scheme structure, and both will be represented by one entry in the scheme. More details on this second type are included in the compendium.

\section{-1 Orthogonal Polynomials descending from\\ the -1 Meixner--Pollaczek polynomials}
In this section, a $q\to-1$ limit of the q-Meixner-Pollaczek family is considered to introduce the $-1$ Meixner-Pollaczek polynomials. These are shown to be a one-parameter extension of
the generalized Hermite polynomials.

\subsection{-1 Meixner-Pollaczek polynomials}

The monic q-Meixner-Pollaczek polynomials are defined by the recurrence relation
\begin{gather}
    x \mathbf{P}_n(x)=\mathbf{P}_{n+1}(x)+a q^{n} \cos \phi \mathbf{P}_n(x)+\frac{1}{4}\left(1-q^{n}\right)\left(1-a^{2} q^{n-1}\right) \mathbf{P}_{n-1}(x),
\end{gather}
with $\mathbf{P}_{-1}(x)=0$ and $\mathbf{P}_0(x)=1$. Multiplying the variable $x$ by $\sqrt{1+q}$, renormalizing the polynomials and using the parametrization 
\begin{gather}
    q=-e^{-\varepsilon},\quad a=-e^{-\varepsilon\pp{\alpha+\frac{1}{2}}},\quad \phi=\frac{\pi}{2}+\sqrt{\varepsilon}\beta,
\end{gather}
we find a set of polynomials that will be denoted $\mathbf{M}_{n}(x;\alpha,\gamma)$ and that will be called the -1 Meixner-Pollaczek. They verify the recurrence relation
\begin{gather}
    x \mathbf{M}_{n}(x)=\mathbf{M}_{n+1}(x)+(-1)^n\gamma \mathbf{M}_{n}(x)+u_n \mathbf{M}_{n-1}(x),\label{Rec_-1MPIn}
\end{gather}
where
\begin{gather}
    u_{2n} = n, \qquad u_{2n+1} = n+\alpha+\frac{1}{2},
\end{gather}
and $\mathbf{M}_{-1}(x)=0$ and $\mathbf{M}_0(x)=1$. Many paths can be used to find the explicit expression of $\mathbf{M}_{n}(x)$. We will obtain it by first noticing (see below) that the recurrence relation of $\mathbf{M}_{n}(x)$ is a limit case of the recurrence relation of the continuous -1 Hahn polynomials. Applying the same limit on the explicit expression of these polynomials will yield the expression of $\mathbf{M}_{n}(x)$. If we define $\hat{\mathbf{K}}_{n}^{\pp{1}}(x;\alpha,\beta,\gamma) = \frac{1}{\pp{2\alpha}^{n/2}}\mathbf{K}_{n}^{\pp{1}}(\sqrt{2\alpha}x;\alpha,\beta,\gamma)$ and parametrize as follows
\begin{gather}
    \alpha = h, \qquad \beta = \sqrt{\frac{h}{2}}b, \qquad \gamma = \frac{2a-1}{4},
\end{gather}
the recurrence relation (\ref{recC-1HIn}) becomes
\begin{gather}
    x \hat{\mathbf{K}}_{n}^{\pp{1}}(x)=\hat{\mathbf{K}}_{n+1}^{\pp{1}}(x)+b_{n}^{\pp{h}} \hat{\mathbf{K}}_{n}^{\pp{1}}(x)+u_{n}^{\pp{h}} \hat{\mathbf{K}}_{n-1}^{\pp{1}}(x),
\end{gather}
where
\begin{gather}
    b_{n}^{\pp{h}}=\left\{\begin{array}{ll}
b\pp{1-\frac{n}{(n+2h+a+1/2)}} & n \text { even } \\
b\pp{1-\frac{n+4h+2a+2}{(n+2h+a+1/2)}} & n \text { odd }
\end{array}\right., \\
    u_{n}^{\pp{h}}=\left\{\begin{array}{ll}
\frac{n(n+4h+2 a+1)(2h)^2}{4(2h)(n+2h+a+1/2)^{2}}\|1+\frac{n+2a+\frac{1}{2}}{2h}+\frac{ib}{\sqrt{2h}}\|^{2} & n \text { even } \\
\frac{(n+4h+1)(n+2a)(2h)^2}{4(2h)(n+2h+a+1/2)^{2}} \|1+\frac{n+2a+\frac{1}{2}}{2h}\|^{2} & n \text { odd }
\end{array}\right..
\end{gather}
The $h\to\infty$ limit is then straightforward to compute, and the result is 
\begin{gather}
    b_{n}^{\pp{\infty}}= (-1)^n b, \hspace{2cm}
    u_{n}^{\pp{\infty}}=\left\{\begin{array}{ll}
\frac{n}{2} & n \text { even } \\[0.2cm]
\frac{n+2a}{2}  & n \text { odd }
\end{array}\right.,
\end{gather}
which coincides with the definition given in (\ref{Rec_-1MPIn}). We thus have established that
\begin{gather}
    \lim\limits_{h\to\infty}\frac{1}{\pp{2h}^\frac{n}{2}}\mathbf{K}_{n}^{(1)}\pp{\sqrt{2h}x,h,\sqrt{\frac{h}{2}}b,\frac{2a-1}{4}}=\mathbf{M}_n\pp{x;a,b}.
\end{gather}
A very similar relation is also valid for the second type of continuous -1 Hahn polynomials. Using this result with the expressions (\ref{CoBIBHR1})-(\ref{CoBIBHR2}) where $\delta$ has been replaced with $\beta$, we find the expressions
\begin{gather}
    \begin{array}{l}
\mathbf{M}_{2n}(x;\alpha,\gamma)=(-1)^{n}\poch{\alpha+1/2}{n} \pFq{1}{1}{-n}{\alpha+1/2}{x^2-\gamma^2}, \\
\mathbf{M}_{2n+1}(x;\alpha,\gamma)=(-1)^{n}\poch{\alpha+3/2}{n} \pp{x-\gamma}\pFq{1}{1}{-n}{\alpha+3/2}{x^2-\gamma^2}.
\end{array}
\end{gather}
These polynomials were already characterized in \cite{CHI} where their orthogonality relation is derived and the second order Dunkl differential operator that they diagonalize is also given. These are included in the compendium. Even though they were already identified as a limit of the Chihara polynomials, the fact that the -1 Meixner-Pollaczek polynomials are, as expected from the given name, a $q\to-1$ limit of the $q$-Meixner-Pollaczek polynomials and that they descendant from the continuous -1 Hahn polynomials is new. 

The natural way to complete this category of $-1$ polynomials is by taking $\gamma=0$ to obtain the generalized Hermite polynomials and then $\alpha = 0$, to obtain the Hermite polynomials. This path of specialization is also coherent with the rest of the scheme, since the limit from the Chihara polynomials to the -1 Meixner-Pollaczek is also valid after the specialization of the Chihara into the generalized Gegenbauer and Gegenbauer polynomials. Then, the generalized Hermite and Hermite polynomials are included as limit cases of the generalized Gegenbauer and Gegenbauer polynomials, respectively.

\section{Continuous Complementary Bannai-Ito and specialization}
In this section, a continuous version of the complementary Bannai-Ito polynomials is investigated with the goal of determining if the Chihara polynomials have a 4-parameters generalization. This family of polynomials is seen to lie beyond the scope of classical orthogonality, since no real three-term recurrence relation can be obtained. Nonetheless, it is shown that some specializations of this family such as the Chihara, the generalized symmetric Bannai-Ito and symmetric Bannai-Ito polynomials possess three-term recurrence relations. The second family is characterized and shown to be related to the Wilson and continuous dual Hahn polynomials.

The continuous Bannai-Ito polynomials were obtained by introducing an imaginary variable and complex parameters in the non-truncated Bannai-Ito polynomials. Using two pairs of complex conjugated parameters, it was shown that the recurrence relation was real and that an infinite family of polynomials orthogonal with respect to a positive continuous measure results from the procedure. The same method will be used here, but with the non-truncated complementary Bannai-Ito polynomials. 

Monic complementary Bannai-Ito polynomials with imaginary variable denoted $\tilde{\mathbf{I}}_{n}(x) = (-i)^n\mathbf{I}_{n}(ix)$, satisfy the recurrence relation
\begin{gather}\label{RecCCBIIn}
    x\tilde{\mathbf{I}}_{n}(x)=\tilde{\mathbf{I}}_{n+1}(x)-i(-1)^{n} \rho_{2} \tilde{\mathbf{I}}_{n}(x)+\tau_{n} \tilde{\mathbf{I}}_{n-1}(x),\\[0.3cm]
    \begin{array}{l}
\tau_{2 n}=\frac{n\left(n+\rho_{1}-r_{1}+1 / 2\right)\left(n+\rho_{1}-r_{2}+1 / 2\right)\left(n-r_{1}-r_{2}\right)}{(2 n+g)(2 n+g+1)}, \\[0.3cm]
\tau_{2 n+1}=\frac{(n+g+1)\left(n+\rho_{1}+\rho_{2}+1\right)\left(n+\rho_{2}-r_{1}+1 / 2\right)\left(n+\rho_{2}-r_{2}+1 / 2\right)}{(2 n+g+1)(2 n+g+2)},
\end{array}
\end{gather}
with $g=\rho_1+\rho_2-r_1-r_2$. It is already clear that $\rho_2$ needs to be imaginary for the recurrence relation to be real. Consider the following parametrization
\begin{gather}\label{changevar}
    \rho_1 = a_1+ib_1-1, \qquad\rho_2 = ib_2, \qquad -r_1 = a_3+ib_3-\frac{1}{2}, \qquad -r_2 = a_4+ib_4-\frac{1}{2}.
\end{gather}
Demanding that the family be positive definite leads to three conditions on the coefficients in the recurrence relation (\ref{RecCCBIIn}):
\begin{enumerate}
  \item $b_1+b_2+b_3+b_4 = 0$,
  \item $\pp{n+a_1+a_3+i\bb{b_1+b_3}}\pp{n+a_1+a_4+i\bb{b_1+b_4}}\pp{n+a_3+a_4-1+i\bb{b_3+b_4}}\in \mathbb{R}$,
  \item $\pp{n+a_3+i\bb{b_2+b_3}}\pp{n+a_4+i\bb{b_2+b_4}}\pp{n+a_1+i\bb{b_1+b_2}}\in \mathbb{R}$.
\end{enumerate}
These three conditions are incompatible if $b_2\neq0$. Solving those equations with $b_2=0$, one finds that $a_3=a_1$, $b_3=-b_1$ and $b_4 = -b_2=0$. We then define the continuous complementary Bannai-Ito polynomials $\tilde{\mathbf{I}}_{n}(x;a_1,b_1,a_2,b_2)$ in this framework as the family of polynomials satisfying the recurrence relation
\begin{gather}\label{RecCCBIInNew}
    x\tilde{\mathbf{I}}_{n}(x)=\tilde{\mathbf{I}}_{n+1}(x)+(-1)^{n} b_{2} \tilde{\mathbf{I}}_{n}(x)+\tau_{n} \tilde{\mathbf{I}}_{n-1}(x),\\[0.3cm]
    \begin{array}{l}
\tau_{2 n}=\frac{n\left(n+2a_1-1\right)\left(n+a_1+a_2-1+i\bb{b_1-b_2}\right)\left(n+a_1+a_2-1-i\bb{b_1+b_2}\right)}{(2 n+2a_1+a_2-2)(2 n+2a_1+a_2-1)}, \\[0.3cm]
\tau_{2 n+1}=\frac{(n+2a_1+a_2-1)\left(n+a_2\right)\left(n+a_1+i\bb{b_1+b_2}\right)\left(n+a_1-i\bb{b_1-b_2}\right)}{(2 n+2a_1+a_2-1)(2 n+2a_1+a_2)},
\end{array}
\end{gather}
with $\tilde{\mathbf{I}}_{-1}(x)=0$ and $\tilde{\mathbf{I}}_{0}(x)=1$. Although these polynomials are not orthogonal with a positive definite measure on the real line since the recurrence coefficients are not real, many specializations and limits of them are orthogonal. They will be included in the scheme with this caveat as they are part of the overall structure, but they shall of course not be put on equal footing with the other orthogonal families. With $l_n$ normalization factors making them monic, they are expressed as follows in terms of the Wilson polynomials $\mathbf{W}_{n}(x;a,b,c,d)$.
\begin{gather}\label{ExpCCBI}
    \begin{array}{l}
\tilde{\mathbf{I}}_{2n}\pp{x;a_1,b_1,a_2,b_2}= l_{2n}\mathbf{W}_{n}(x^2;ib_2,a_1+ib_1,a_1-ib_1,a_2-ib_2), \\
\tilde{\mathbf{I}}_{2n+1}(x;a_1,b_1,a_2,b_2)= l_{2n+1}(x-b_2)\mathbf{W}_{n}(x^2;1+ib_2,a_1+ib_1,a_1-ib_1,a_2-ib_2),
\end{array}
\end{gather}
where 
\begin{align}
    l_{2n} &= \frac{\pp{-1}^n\poch{a_1+i\bb{b_2+b_1}}{n}\poch{a_1+i\bb{b_2-b_1}}{n}\poch{a_2}{n}}{\poch{n+2a_1+a_2-1}{n}},\\ l_{2n+1} &= \frac{\pp{-1}^n\poch{1+a_1+i\bb{b_2+b_1}}{n}\poch{1+a_1+i\bb{b_2-b_1}}{n}\poch{1+a_2}{n}}{\poch{n+2a_1+a_2}{n}}.
\end{align}

\subsection{A Limit to the Chihara Polynomials}
The Chihara polynomials can be obtained as a limit of the continuous complementary Bannai-Ito polynomials. Let $b_1$ and $b_2$ be parametrized as follows
\begin{gather}
    b_1 = hc_1, \qquad b_2 = hc_2,
\end{gather}
and denote the renormalized polynomials by
\begin{gather}
    \tilde{\mathbf{I}}_{n}^{(h)}(x) = h^{-n}\pp{c_1^2-c_2^2}^{\frac{-n}{2}} \tilde{\mathbf{I}}_{n}\pp{h\sqrt{c_1^2-c_2^2}x}.
\end{gather}
Upon taking the limit $h\to\infty$, the recurrence relation (\ref{RecCCBIInNew}) goes to
\begin{gather}\label{NewRecC}
    x\tilde{\mathbf{I}}_{n}^{(\infty)}(x)=\tilde{\mathbf{I}}_{n+1}^{(\infty)}(x)+(-1)^{n} \pp{\frac{c_2}{\sqrt{c_1^2-c_2^2}}} \tilde{\mathbf{I}}_{n}^{(\infty)}(x)+u_{n} \tilde{\mathbf{I}}_{n-1}^{(\infty)}(x),\\[0.3cm]
    \begin{array}{l}
u_{2 n}=\frac{n\left(n+2a_1-1\right)}{(2 n+2a_1+a_2-2)(2 n+2a_1+a_2-1)},\qquad 
u_{2 n+1}=\frac{(n+2a_1+a_2-1)\left(n+a_2\right)}{(2 n+2a_1+a_2-1)(2 n+2a_1+a_2)}
\end{array},
\end{gather}
which is the recurrence relation of the Chihara polynomials allowing for an affine transformation of the parameters. We conclude that the Chihara polynomials are limits of the continuous complementary Bannai-Ito polynomials. This corresponds to
\begin{gather}\label{LimCCBItoC}
    \lim\limits_{h\to\infty}\frac{1}{h^n\pp{c_1^2-c_2^2}^{\frac{n}{2}}}\tilde{\mathbf{I}}_{n}\pp{h\sqrt{c_1^2-c_2^2}x;\frac{\beta+1}{2},hc_1,\alpha+1,hc_2}=\mathbf{C}_n\pp{x;\alpha,\beta,\frac{c_2}{\sqrt{c_1^2-c_2^2}}}.
\end{gather}

\subsection{Specialization to Generalized Symmetric Bannai-Ito}\label{secunnammed1}
As mentioned before, taking $b_2=0$ in (\ref{RecCCBIInNew}) gives a real recurrence relation for polynomials that are orthogonal with respect to a positive continuous measure when $\tau_n>0$. These polynomials, that we shall call generalized symmetric Bannai-Ito polynomials and denote $\hat{\mathbf{I}}_{n}(x)$, are clearly symmetric. They are a made out of Wilson polynomials, and thus the Chihara method for symmetric moment functional \cite{CHIbook} can be used on the expressions (\ref{ExpCCBI}) with $b_2=0$. For the remainder of the section, we will use the change of parameter $a=a_1+ib_1$, $b=a_1-ib_1$, $c=a_2$ and denote the polynomials by $\hat{\mathbf{I}}_{n}(x;a,b,c)$. Applying the Chihara procedure, the weight function for the generalized symmetric Bannai-Ito family is found to be
\begin{gather}
    \omega(x) = \left|\frac{\Gamma\pp{ix}\Gamma\pp{a+ix}\Gamma\pp{b+ix}\Gamma\pp{c+ix}}{\Gamma\pp{2ix}}\right|^2,
\end{gather}
on the interval $[-\infty,\infty]$. When Re$\pp{a,b,c}>0$, $a+b+c>1$ and non-real parameters occur in conjugate pairs, the orthogonality relation reads
\begin{gather}\label{OrtUN1IN}
    \frac{1}{4\pi}\int_{-\infty}^{\infty} \omega(x) \hat{\mathbf{I}}_{n}(x)\hat{\mathbf{I}}_{m}(x) \mathrm{d} x=\kappa_n \delta_{n m}, \\[0.3cm]
    \kappa_{n} = \frac{\Gamma\pp{n+a+b}\Gamma\pp{n+a+c}\Gamma\pp{n+b+c}\Gamma\pp{n+a}\Gamma\pp{n+b}\Gamma\pp{n+c}n!}{\Gamma\pp{2n+a+b+c}\poch{n+a+b+c-1}{n}}.
\end{gather}

The generalized symmetric Bannai-Ito polynomials are eigenfunctions of a modified version of the difference operator of the complementary Bannai-Ito polynomials \cite{CBI} where the shift is on the imaginary axis $S^{\pm}f(x) = f(x\pm i)$ in view of the change of variable introduced in (\ref{RecCCBIIn}) and where the parameters are taken according to (\ref{changevar}): 
\begin{gather}
    D_\sigma \hat{\mathbf{I}}_{n}(x) = \Lambda_{n}^{(\sigma)}\hat{\mathbf{I}}_{n}(x)\label{D0diffop},\\[0.3cm]
    \begin{array}{l}
        \Lambda_{2n}^{(\sigma)} = n^2 +\pp{a+b+c-1}n, \qquad
        \Lambda_{2n+1}^{(\sigma)} = n^2 +\pp{a+b+c}n+\sigma,
    \end{array}\\[0.3cm]
    D_{\sigma} = D_0 + \frac{\sigma}{2}\pp{I-R},\\[0.3cm]
    D_0 = B(x)S^++A(x)S^-+C(x)R-\pp{A(x)+B(x)+C(x)}I,
\end{gather}
\begin{gather}
    A(x) = \frac{\pp{ix+a}\pp{ix+b}\pp{ix+c}}{2\pp{2ix+1}},\\[0.3cm]
    B(x) = \frac{\pp{ix-a}\pp{ix-b}\pp{ix-c}}{2\pp{2ix-1}},\\[0.3cm]
    C(x) = \frac{1}{2}\pp{ab+ac+bc-x^2}-A(x)-B(x).
\end{gather}

This establishes the bispectrality property of the generalized symmetric Bannai-Ito and complete their characterization. It is also interesting to note that the limit (\ref{LimCCBItoC}) is also valid for $\hat{\mathbf{I}}_{n}(x)$ with $c_2=0$, which relates the generalized symmetric Bannai-Ito family to the generalized Gegenbauer polynomials.

\subsection{A limit to Symmetric Bannai-Ito polynomials}
In view of \cref{secunnammed1}, the limit from the Wilson to the continuous dual Hahn polynomials can be used similarly for their $-1$ counterpart. This amount to taking the limit when $c$ goes to $\infty$. One then obtains a family of orthogonal polynomials called symmetric Bannai-Ito with complex parameters and denoted $\hat{\mathbf{S}}_{n}(x)$. They were introduced in a different but equivalent way in \cite{BI}. They are given by 
\begin{gather}\label{ExpS}
    \begin{array}{l}
\hat{\mathbf{S}}_{2n}\pp{x;a,b}= k_{2n}\mathbf{S}_{n}(x^2;0,a,b), \\
\hat{\mathbf{S}}_{2n+1}(x;a,b)= k_{2n+1}x\mathbf{S}_{n}(x^2;1,a,b),
\end{array}
\end{gather}
with normalization constants $k_n$ ensuring that they are monic and where $\mathbf{S}_{n}(x)$ are continuous dual Hahn polynomials.
These polynomials satisfy the recurrence relation
\begin{gather}\label{RecSInNew}
    x\hat{\mathbf{S}}_{n}(x)=\hat{\mathbf{S}}_{n+1}(x)+\tau_{n} \hat{\mathbf{S}}_{n-1}(x),\\[0.3cm]
    \begin{array}{l}
\tau_{2 n}=n\left(n+a+b-1\right),\qquad
\tau_{2 n+1}=\pp{n+a}\pp{n+b}.
\end{array}
\end{gather}
When Re$\pp{a,b}>0$ and non-real parameters occur in conjugated pairs, the orthogonality relation reads
\begin{gather}\label{OrtUN2IN}
    \frac{1}{4\pi}\int_{-\infty}^{\infty} \omega(x) \hat{\mathbf{S}}_{n}(x)\hat{\mathbf{S}}_{m}(x) \mathrm{d} x=\kappa_n \delta_{n m} ,\\[0.3cm]
    \kappa_{n} = \Gamma\pp{n+a+b}\Gamma\pp{n+a}\Gamma\pp{n+b}n!,
\end{gather}
with 
\begin{gather}
    \omega(x) = \left|\frac{\Gamma\pp{ix}\Gamma\pp{a+ix}\Gamma\pp{b+ix}}{\Gamma\pp{2ix}}\right|^2.
\end{gather}
The difference operator (\ref{D0diffop}) exists in the limit $c\to\infty$ if we divide on both side by $c$ before taking the limit. The following difference equation is obtained for $\hat{\mathbf{S}}_n\pp{x}$. 
\begin{gather}
    D_\sigma \hat{\mathbf{S}}_{n}(x) = \Lambda_{n}^{(\sigma)}\hat{\mathbf{S}}_{n}(x),\label{D0diffopS}\\[0.3cm]
    \begin{array}{l}
        \Lambda_{2n}^{(\sigma)} = n ,\qquad
        \Lambda_{2n+1}^{(\sigma)} = n+\sigma,
    \end{array}\\[0.3cm]
    D_{\sigma} = D_0 + \frac{\sigma}{2}\pp{I-R},\\[0.3cm]
    D_0 = B(x)S^++A(x)S^-+C(x)R-\pp{A(x)+B(x)+C(x)}I,
\end{gather}
\begin{gather}
    A(x) = \frac{\pp{ix+a}\pp{ix+b}}{2\pp{1+2ix}},\\[0.3cm]
    B(x) = \frac{\pp{ix-a}\pp{ix-b}}{2\pp{1-2ix}},\\[0.3cm]
    C(x) = \frac{a+b}{2}-A(x)-B(x).
\end{gather}
We obtained the symmetric Bannai-Ito by first taking $b_2=0$ and then $c\to\infty$. If we reverse the order of these operations on the polynomials, we can observe that the polynomials $\hat{\mathbf{S}}_n$ are some specialization of a $q\to-1$ limit of the continuous dual q-Hahn polynomials. This observation is due to the correspondence between the $c\to\infty$ limit and the limit relation between Askey-Wilson and continuous dual q-Hahn polynomials. It should be noted that the symmetric Bannai-Ito polynomials can also be obtained by setting $\beta=0$ in the continuous $-1$ Hahn polynomials of \cref{C-1HP}. Finally, the symmetric Bannai-Ito have the generalized Hermite polynomials as a limit. This connection is found in the compendium.

\section{Conclusion}

In this paper, we have constructed the continuous part of the -1 orthogonal polynomials scheme and presented it using four categories. First, we completed the bulk of the scheme, which contains the -1 polynomials obtained from the big $q$-Jacobi polynomials. The continuous Bannai-Ito polynomials were identified as the top family of the continuous part of the scheme and as parent of the big -1 Jacobi polynomials. Furthermore, they were identified as a generalization of the continuous -1 Hahn polynomials. These polynomials and the continuous Bannai-Ito polynomials formed the second category. In the third part, two new ways of looking at the -1 Meixner-Pollaczek were proposed. They are seen on the one hand as a $q\to-1$ limit of the q-Meixner-Pollaczek polynomials, and on the other as a limit of the novel continuous $-1$ Hahn polynomials. The $-1$ Meixner-Pollaczek polynomials, together with the generalized Hermite and Hermite polynomials, made out the third category of $-1$ continuous polynomials. In the later part, a continuous equivalent to the complementary Bannai-Ito was looked at. Even if the resulting polynomials are not orthogonal, many specializations and limits are. It was found that the Chihara polynomials are descendants of this a family. Another new family was introduced: the generalized symmetric Bannai-Ito polynomials, that are also descendants of the continuous complementary Bannai-Ito polynomials. The explicit expression, orthogonality relation and difference equation were obtained using the properties of the Wilson polynomials. A compendium of the main properties of all the families included in the continuous part of the scheme is presented in the appendix. To arrive at a complete $-1$ scheme, there remains to determine its discrete part. While many elements of it have been identified, the full picture requires more work. It is our plan to supplement the continuous part of the $-1$ scheme presented here with its discrete complement in an upcoming report.

\subsubsection*{Acknowledgements}
JP holds a scholarship from Fonds de recherche québécois – nature et technologies (FRQNT) and an academic excellence scholarships from Hydro-Québec. The research of LV is supported in part by a Discovery grant from the Natural Science and Engineering Research Council (NSERC) of Canada.

\newpage
\newgeometry{bottom=2cm,top=2cm,left=1cm,right=1.5cm}
\begin{center}
\begin{Huge}
Scheme\\
of\\
Continuous -1 Hypergeometric\\
Orthogonal Polynomials\vspace{1cm}
\end{Huge}
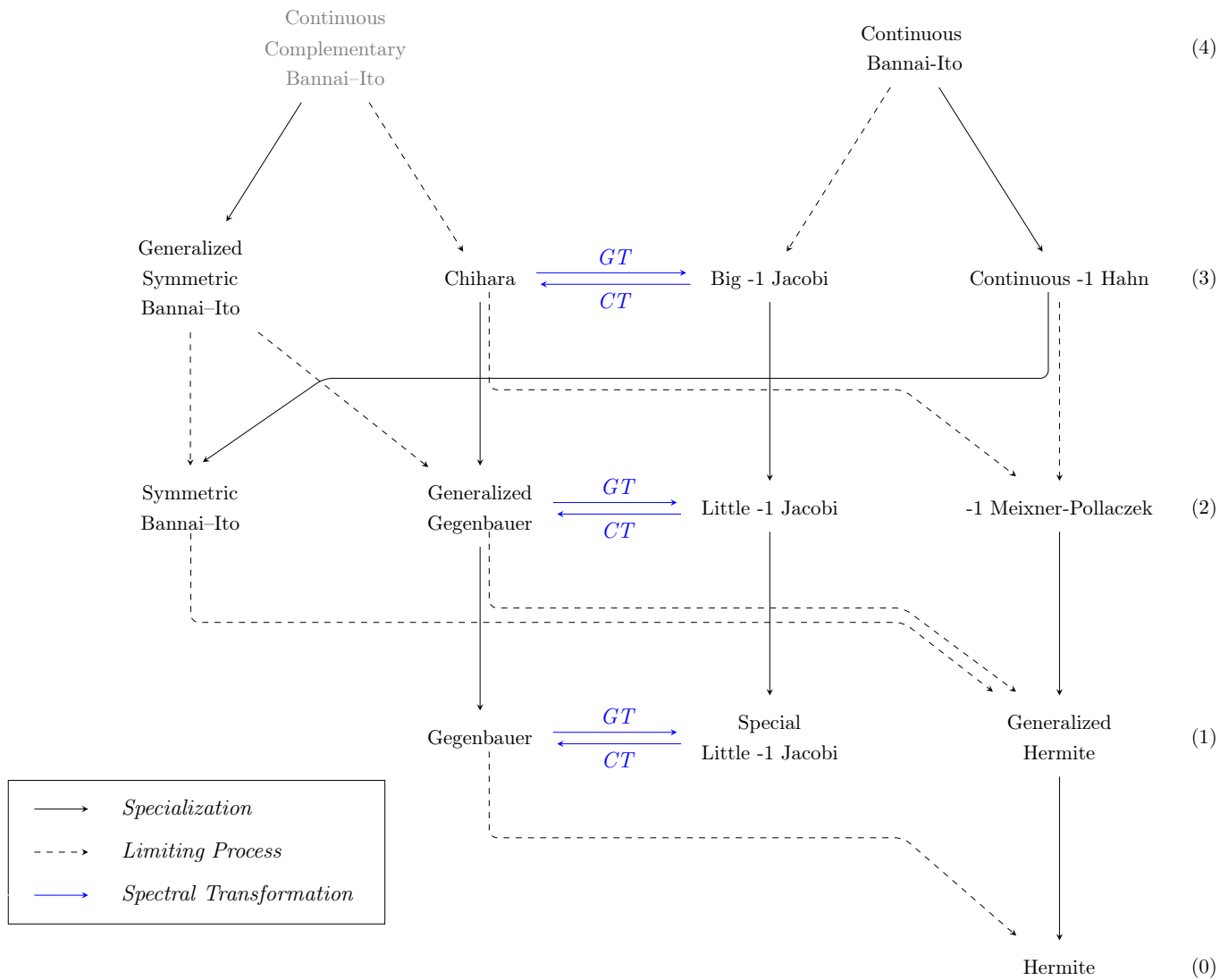
\begin{figure}[H]
\begin{adjustbox}{width=\textwidth}
\begin{tikzpicture}[>=stealth,->,shorten >=2pt,looseness=.5,auto]
\tikzstyle{every node}=[font=\footnotesize]
\matrix [matrix of math nodes,
column sep={2.5cm,between origins},
row sep={4cm,between origins}]
{
|(11)| & |(12)| \begin{tabular}{c} \( \text{\color{gray}Continuous\color{black}} \) \\ \( \text{\color{gray}Complementary\color{black}} \) \\ \( \text{\color{gray}Bannai–Ito\color{black}}\) \end{tabular} & |(13)| &|(14)|& |(15)| & |(16)|  \begin{tabular}{c} \( \text{Continuous } \) \\ \( \text{Bannai-Ito} \) \end{tabular} &|(17)|& |(18)|  (4)\\
|(21)| \begin{tabular}{c} \( \text{Generalized} \) \\ \( \text{Symmetric} \) \\ \( \text{Bannai–Ito}\) \end{tabular} & |(22)| & |(23)|  \begin{tabular}{c} \( \text{Chihara} \) \end{tabular} & |(24)| & |(25)| \begin{tabular}{c} \( \text{Big -1 Jacobi} \)\end{tabular}  &|(26)| &|(27)| \begin{tabular}{c} \( \text{Continuous -1 Hahn} \) \end{tabular} & |(28)| (3) \\
|(31)| \begin{tabular}{c}  \( \text{Symmetric} \) \\ \( \text{Bannai–Ito}\) \end{tabular} & |(32)| & |(33)|  \begin{tabular}{c} \( \text{Generalized} \) \\\( \text{Gegenbauer}\) \end{tabular} & |(34)| & |(35)| \begin{tabular}{c} \( \text{Little -1 Jacobi} \)  \end{tabular}  &|(36)| &|(37)| \begin{tabular}{c} \( \text{-1 Meixner-Pollaczek} \) \end{tabular} & |(38)| (2) \\
|(41)|  & |(42)| & |(43)|  \begin{tabular}{c} \( \text{Gegenbauer}\) \end{tabular} & |(44)| & |(45)| \begin{tabular}{c} \(\text{Special}\)\\\( \text{Little -1 Jacobi} \)\end{tabular}  &|(46)| &|(47)|  \begin{tabular}{c} \( \text{Generalized} \) \\\( \text{Hermite}\) \end{tabular}  & |(48)| (1) \\
|(51)|  & |(52)| & |(53)|   & |(54)| & |(55)|  &|(56)| &|(57)|  \begin{tabular}{c}\( \text{Hermite}\) \end{tabular}  & |(58)| (0) \\
};
\tikzstyle{every node}=[font=\small\itshape]
\draw[blue] ([shift={(0,0.1)}]23.east) -- ([shift={(0,0.1)}]25.west) node[midway,above] {GT};
\draw[blue] ([shift={(0,-0.1)}]25.west) -- ([shift={(0,-0.1)}]23.east) node[midway,below] {CT};
\draw[blue] ([shift={(0,0.1)}]33.east) -- ([shift={(0,0.1)}]35.west) node[midway,above] {GT};
\draw[blue] ([shift={(0,-0.1)}]35.west) -- ([shift={(0,-0.1)}]33.east) node[midway,below] {CT};
\draw[blue] ([shift={(0,0.1)}]43.east) -- ([shift={(0,0.1)}]45.west) node[midway,above] {GT};
\draw[blue] ([shift={(0,-0.1)}]45.west) -- ([shift={(0,-0.1)}]43.east) node[midway,below] {CT};
\draw[dashed] (12) -- (23);
\draw[dashed] (16) -- (25);
\draw (12) -- (21);
\draw (16) -- (27);
\draw[dashed] (21) -- (31);
\draw (23) -- (33);
\draw (25) -- (35);
\draw[dashed] (27) -- (37);
\draw[dashed] (21) -- (33);
\draw (33) -- (43);
\draw (35) -- (45);
\draw (37) -- (47);
%
\draw (47) -- (57);
\draw[rounded corners,dashed] (-3.15,3.5) -- (-3.15,1.8)--(4.0,1.8)--(6,0.25);
\draw[rounded corners,dashed] (-3.15,-0.67) -- (-3.15,-2)--(4,-2)--(6,-3.50);
\draw[rounded corners,dashed] (-3.15,-4.5) -- (-3.15,-6)--(4,-6)--(6,-7.75);
\draw[rounded corners] (6.5,3.5) -- (6.5,2) -- (-6,2) -- (-8.15,0.5);
\draw[rounded corners,dashed] (-8.3,-0.7) -- (-8.3,-2.25)--(4,-2.25)--(5.6,-3.5);
\draw[->]  (11-22,-5.5) -- (12-22,-5.5);
\node[align=left,text width=5cm] at (15-22,-5.5) {Specialization};
\draw[->,dashed]  (11-22,-6.25) -- (12-22,-6.25);
\node[align=left,text width=5cm] at (15-22,-6.25) {Limiting Process};
\draw[->,blue]  (11-22,-7) -- (12-22,-7);
\node[align=left,text width=5cm] at (15-22,-7) {Spectral Transformation};
\node[rectangle,draw,minimum width = 6.5cm,minimum height = 2.5cm] (r) at (13.8-22,-6.25) {};
\draw[->,white]  (11-23,-7) -- (12-23,-7);
\end{tikzpicture}
\end{adjustbox}
\end{figure}

\thispagestyle{empty}
\end{center}
\restoregeometry

\textwidth 170mm 
\textheight 225mm 
\topmargin -18mm
\oddsidemargin -0.5cm

\appendix
\section{A compendium of the properties of the\\ \texorpdfstring{$-1$}{-1} continuous orthogonal polynomials}
\subsection{Continuous Bannai-Ito}
\textbf{\textit{\large Hypergeometric Representation}}\\
\begin{flalign}
    &\frac{\mathbf{Q}_{2n}\left(x;\alpha,\beta,\gamma,\delta\right)}{\pp{-2i}^{2n}} = &\end{flalign}\vspace{-0.75cm}\begin{align}\nonumber
    &\xi_{2n}\kappa_{n-1}^{(2)}\left(\frac{ix}{2}-\mathfrak{b}-\frac{1}{2}\right)\pFq{4}{3}{-n+1,n+g+2,\frac{ix}{2}+\mathfrak{b}+1,\frac{-ix}{2}+\mathfrak{b}+\frac{3}{2}}{\frac{5}{2}+\mathfrak{a}+\mathfrak{b},2+\mathfrak{b}+\mathfrak{c},2+\mathfrak{b}+\mathfrak{d}}{1} + \\
    &\kappa_n^{(1)}\pFq{4}{3}{-n,n+g+1,\frac{ix}{2}+\mathfrak{b},\frac{-ix}{2}+\mathfrak{b}+\frac{1}{2}}{\frac{3}{2}+\mathfrak{a}+\mathfrak{b},1+\mathfrak{b}+\mathfrak{c},1+\mathfrak{b}+\mathfrak{d}}{1}, \nonumber\end{align}
    \begin{flalign} &\frac{\mathbf{Q}_{2n+1}\left(x;\alpha,\beta,\gamma,\delta\right)}{\pp{-2i}^{2n+1}} =&\end{flalign}\vspace{-0.75cm}\begin{align}\nonumber &\kappa_{n}^{(2)}\left(\frac{ix}{2}-\mathfrak{b}-\frac{1}{2}\right)\pFq{4}{3}{-n,n+g+2,\frac{ix}{2}+\mathfrak{b}+1,\frac{-ix}{2}+\mathfrak{b}+\frac{3}{2}}{\frac{5}{2}+\mathfrak{a}+\mathfrak{b},2+\mathfrak{b}+\mathfrak{c},2+\mathfrak{b}+\mathfrak{d}}{1}+\\
    &\eta_{2n+1}\kappa_n^{(1)}\pFq{4}{3}{-n,n+g+1,\frac{ix}{2}+\mathfrak{b},\frac{-ix}{2}+\mathfrak{b}+\frac{1}{2}}{\frac{3}{2}+\mathfrak{a}+\mathfrak{b},1+\mathfrak{b}+\mathfrak{c},1+\mathfrak{b}+\mathfrak{d}}{1}, \nonumber
\end{align}
where
\begin{gather}
    g=2\alpha+2\gamma+1,\qquad \mathfrak{a}=\overline{\mathfrak{d}}=\alpha+i\beta, \qquad \mathfrak{b}=\overline{\mathfrak{c}}=\gamma+i\delta,
\end{gather}
\begin{gather}
\begin{array}{ll}
\xi_{2n}=\frac{n\pp{n+\mathfrak{c}+\mathfrak{d}+\frac{1}{2}}}{\left(2n+g\right)}, & \kappa_{n}^{(1)}=\frac{\poch{\frac{3}{2}+\mathfrak{a}+\mathfrak{b}}{n} \poch{1+\mathfrak{b}+\mathfrak{c}}{n} \poch{1+\mathfrak{b}+\mathfrak{d}}{n}}{\poch{n+g+1}{n}}, \\[0.15cm]
\eta_{2n+1}=\frac{\pp{n+\mathfrak{b}+\mathfrak{c}+1}\pp{n+\mathfrak{b}+\mathfrak{d}+1}}{\pp{2n+g+1}}, &  \kappa_{n}^{(2)}=\frac{\poch{\frac{5}{2}+\mathfrak{a}+\mathfrak{b}}{n} \poch{2+\mathfrak{b}+\mathfrak{c}}{n} \poch{2+\mathfrak{b}+\mathfrak{d}}{n}}{\poch{n+g+2}{n}}.
\end{array}
\end{gather}
\textbf{\textit{\large Orthogonality Relation}}\\
If $\alpha,\beta,\gamma,\delta\in\mathbb{R}^+$
\begin{gather}
    \int_{-\infty}^{\infty} W(x) \mathbf{Q}_{n}(x) \mathbf{Q}_{m}(x) \mathrm{d} x=4\pi h_{0}\kappa_n \delta_{n m}, \\[0.3cm]
    W(x)=\left|\frac{\Gamma(\mathfrak{a}+i x / 2+1) \Gamma(\mathfrak{b}+i x / 2+1) \Gamma(\mathfrak{c}+i x / 2+1 / 2) \Gamma(\mathfrak{d}+i x / 2+1 / 2)}{\Gamma(1 / 2+i x)}\right|^2,\\[0.3cm]
    h_{0}=\frac{\Gamma(\mathfrak{a}+\mathfrak{b}+3 / 2) \Gamma(\mathfrak{a}+\mathfrak{c}+1) \Gamma(\mathfrak{b}+\mathfrak{c}+1) \Gamma(\mathfrak{a}+\mathfrak{d}+1) \Gamma(\mathfrak{b}+\mathfrak{d}+1) \Gamma(\mathfrak{c}+\mathfrak{d}+3 / 2)}{\Gamma(\mathfrak{a}+\mathfrak{b}+\mathfrak{c}+\mathfrak{d}+2)},\end{gather}
\begin{minipage}{\textwidth}
    \begin{flalign}
    &\kappa_{2n} = \frac{4^{2n}\Gamma\pp{n+1}\poch{2\alpha+1}{n}\poch{2\gamma+1}{n}}{\poch{2\alpha+2\gamma+2}{2n} \poch{n+2\alpha+2\gamma+2}{n}} \times&\end{flalign}\vspace{-0.70cm}\begin{align}& \bb{\prod\limits_{k=1}^n \| k+\alpha+\gamma+i\pp{\beta-\delta}\|}^2 \bb{\prod\limits_{k=1}^n \| k+\alpha+\gamma+\frac{1}{2}+i\pp{\beta+\delta}\|}^2,\nonumber\end{align}
    \end{minipage}
    \begin{minipage}{\textwidth}
    \begin{flalign}
    \kappa_{2n+1} =& \frac{4^{2n+1}\Gamma\pp{n+1}\poch{2\alpha+1}{n+1}\poch{2\gamma+1}{n+1}}{\poch{2\alpha+2\gamma+2}{2n+1} \poch{n+2\alpha+2\gamma+2}{n+1}}\times& \end{flalign}\vspace{-0.70cm}\begin{align} &\bb{\prod\limits_{k=1}^{n+1} \| k+\alpha+\gamma+i\pp{\beta-\delta}\|}^2 \bb{\prod\limits_{k=1}^n \| k+\alpha+\gamma+\frac{1}{2}+i\pp{\beta+\delta}\|}^2.\nonumber
\end{align}\end{minipage}\\[0.5cm]
\textbf{\textit{\large Normalized Recurrence Relation}}\\
\begin{gather}
    x\mathbf{Q}_n(x) = \mathbf{Q}_{n+1}(x) + b_n \mathbf{Q}_n(x) + u_n \mathbf{Q}_{n-1}(x),
\end{gather}
\begin{gather}
    b_{n}=\left\{\begin{array}{ll}
2 \beta-\frac{(n+4 \alpha+2)(\beta-\delta)}{(n+2 \alpha+2 \gamma+2)}-\frac{n(\beta+\delta)}{(n+2 \alpha+2 \gamma+1)} & n \text { even } \\[0.3cm]
2 \beta-\frac{(n+4 \alpha+4 \gamma+3)(\beta+\delta)}{n+2 \alpha+2 \gamma+2}-\frac{(n+4 \gamma+1)(\beta-\delta)}{n+2 \alpha+2 \gamma+1} & n \text { odd }
\end{array}\right.,\\[0.3cm]
u_{n}=\left\{\begin{array}{ll}
\frac{n(n+4 \alpha+4 \gamma+2)\|n+2[\alpha+\gamma+i(\beta+\delta)]+1\|^{2}}{4(n+2 \alpha+2 \gamma+1)^{2}} & n \text { even } \\[0.3cm]
\frac{(n+4 \alpha+1)(n+4 \gamma+1)\|n+2[\alpha+\gamma+i(\beta-\delta)]+1\|^{2}}{4(n+2\alpha+2 \gamma+1)^{2}} & n \text{ odd}
\end{array}\right..
\end{gather}\\
\textbf{\textit{\large Difference Equation}}\\
\begin{gather}
    L \mathbf{Q}_{n}(x)=\lambda_{n} \mathbf{Q}_{n}(x), \quad \lambda_{n}=(-1)^{n}(n+2\alpha+2\gamma+3 / 2),\\[0.3cm]
L =\mathcal{A}\left(S^{+} R-I\right)
+\overline{\mathcal{A}}\left(S^{-} R-I\right)+(2\alpha+2\gamma+3 / 2) I,\\[0.3cm]
\mathcal{A} = \frac{\pp{2 \alpha+1+i\bb{\beta- x}}\pp{2 \gamma+1+i\bb{\delta- x}}}{1-2 i x},
\end{gather}
where $\overline{\mathcal{A}}$ is the complex conjugate of $\mathcal{A}$, $S^{\pm}f(x)=f(x\pm i)$ and $Rf(x)=f(-x)$.\\[0.5cm]
\begin{minipage}{\textwidth}
\textbf{\textit{\large Limit Relations}}\\[0.3cm]
\textbf{\textit{\large Continuous Bannai-Ito $\to$ Big -1 Jacobi}}\\[0.3cm]
The Big -1 Jacobi polynomial can be obtained from the Continuous Bannai-Ito polynomials using the parametrization $\beta\to\frac{\beta}{h}$, $\delta\to\frac{\delta}{h}$, the scaling $x\to2\beta x$ and taking the limit $h\to0$:
\begin{gather}
    \lim\limits_{h\to0}\left(\frac{h}{2\beta}\right)^n \mathbf{Q}_n\left(\frac{2\beta x}{h};\alpha,\frac{\beta}{h},\gamma,\frac{\delta}{h}\right) = \mathbf{J}_n\left(x;4\alpha+1,4\gamma+1,\frac{-\delta}{\beta}\right).
\end{gather}
\end{minipage}\\[0.5cm]
\textbf{\textit{\large Continuous Bannai-Ito $\to$ First type Continuous -1 Hahn}}\\[0.3cm]
The first type of Continuous -1 Hahn polynomials can be obtained from the Continuous Bannai-Ito polynomials by a specialization:
\begin{gather}
    \mathbf{Q}_n\left(x;\alpha,\beta,\gamma,\beta\right) = \mathbf{K}_n^{(1)}\left(x;\alpha,\beta,\gamma\right).
\end{gather}\\
\textbf{\textit{\large Continuous Bannai-Ito $\to$ Second type Continuous -1 Hahn}}\\[0.3cm]
The second type of Continuous -1 Hahn polynomials can be obtained from the Continuous Bannai-Ito polynomials by a specialization:
\begin{gather}
    \mathbf{Q}_n\left(x;\alpha,\beta,\gamma,-\beta\right) = \mathbf{K}_n^{(2)}\left(x;\alpha,\beta,\gamma\right).
\end{gather}

\subsection{Big -1 Jacobi}
\textbf{\textit{\large Hypergeometric Representation}}\\
\begin{flalign}
    &\frac{\mathbf{J}_{2n}\left(x;\alpha,\beta,c\right)}{\eta_{2n}}=&\label{B-1JBHR13}\end{flalign}\vspace{-0.75cm}\begin{align}\nonumber &\pFq{2}{1}{-n,\frac{2n+\alpha+\beta+2}{2}}{\frac{\alpha+1}{2}}{\frac{1-x^2}{1-c^2}}+ 
    \frac{2n(1-x)}{(1+c)(1+\alpha)}\pFq{2}{1}{1-n,\frac{2n+\alpha+\beta+2}{2}}{\frac{\alpha+3}{2}}{\frac{1-x^2}{1-c^2}},  \end{align}
    \begin{flalign}
    &\frac{\mathbf{J}_{2n+1}\left(x;\alpha,\beta,c\right)}{\eta_{2n+1}}=& \label{B-1JBHR23}\end{flalign} \begin{align}\nonumber &\pFq{2}{1}{-n,\frac{2n+\alpha+\beta+2}{2}}{\frac{\alpha+1}{2}}{\frac{1-x^2}{1-c^2}}-
    \frac{(2n+\alpha+\beta+2)(1-x)}{(1+c)(1+\alpha)}\pFq{2}{1}{-n,\frac{2n+\alpha+\beta+4}{2}}{\frac{\alpha+3}{2}}{\frac{1-x^2}{1-c^2}},
\end{align}
where
\begin{gather}
    \eta_{2n} =  \frac{\left(1-c^2\right)^n\poch{\frac{\alpha+1}{2}}{n}}{\poch{\frac{2n+\alpha+\beta+2}{2}}{n}},\qquad \eta_{2n+1} = \frac{\left(1+c\right)\left(1-c^2\right)^{n}\poch{\frac{\alpha+1}{2}}{n+1}}{\poch{\frac{2n+\alpha+\beta+2}{2}}{n+1}}.
\end{gather}
\textbf{\textit{\large Orthogonality Relation}}\\
If $\alpha>0$, $\beta>0$ and $0\leq c<1$
\begin{gather}
     \int_\mathcal{C} \omega(x) \mathbf{J}_n(x)\mathbf{J}_m(x) \mathrm{d} x=\frac{\pp{1-c}\pp{1-c^2}^{\frac{\alpha+\beta}{2}}\Gamma\pp{\frac{\alpha+1}{2}}\Gamma\pp{\frac{\beta+1}{2}}}{\Gamma\pp{\frac{\alpha}{2}+\frac{\beta}{2}+1}} \kappa_n \delta_{n m}, \\[0.3cm]
     \mathcal{C} = \bb{-1,-c}\bigcup\bb{c,1},\\[0.3cm]
     \omega(x) = \theta(x)\frac{\pp{1+x}}{\pp{c+x}}\pp{1-x^2}^{\frac{\alpha-1}{2}}\pp{x^2-c^2}^{\frac{\beta+1}{2}},\\[0.3cm]
    \begin{array}{l}
    \kappa_{2n} = \frac{\Gamma\pp{n+1}\poch{\frac{\alpha+1}{2}}{n}\poch{\frac{\beta+1}{2}}{n}}{\poch{1+\frac{\alpha+\beta}{2}}{2n}\poch{n+1+\frac{\alpha+\beta}{2}}{n}} \pp{1-c^2}^n, \\[0.3cm]
    \kappa_{2n+1} = \frac{\Gamma\pp{n+1}\poch{\frac{\alpha+1}{2}}{n+1}\poch{\frac{\beta+1}{2}}{n+1}}{\poch{1+\frac{\alpha+\beta}{2}}{2n+1}\poch{n+1+\frac{\alpha+\beta}{2}}{n+1}} \pp{1+c}^2\pp{1-c^2}^n .
\end{array}
\end{gather}\\
\textbf{\textit{\large Normalized Recurrence Relation}}\\
\begin{gather}
    x\mathbf{J}_n(x) = \mathbf{J}_{n+1}(x) + \left(1-A_n-C_n\right)\mathbf{J}_n(x) + A_{n-1}C_n \mathbf{J}_{n-1}(x),
\end{gather}
\begin{minipage}{.45\textwidth}
        \begin{eqnarray}
    A_{n}=\left\{\begin{array}{ll}
\frac{\pp{1+c}\pp{n+\alpha+1}}{2n+\alpha+\beta+2}&n\text{ even}\\[0.15cm]
\frac{\pp{1-c}\pp{n+\alpha+\beta+1}}{2n+\alpha+\beta+2}& n\text{ odd}
\end{array}\right.,\nonumber
\end{eqnarray}
    \end{minipage}%
    \begin{minipage}{0.55\textwidth}
       \begin{eqnarray}
   C_{n}=\left\{\begin{array}{ll}
\frac{\pp{1-c}n}{2n+\alpha+\beta}&n\text{ even}\\[0.15cm]
\frac{\pp{1+c}\pp{n+\beta}}{2n+\alpha+\beta}& n\text{ odd}
\end{array}\right..\label{B-1JRRC123}
\end{eqnarray}
    \end{minipage}\\[0.5cm]
\textbf{\textit{\large Differential Equation}}\\
\begin{gather}
    L \mathbf{J}_{n}(x)=\lambda_{n} \mathbf{J}_{n}(x), \quad \lambda_{n}=\left\{\begin{array}{ll}
-2n&n\text{ even}\\[0.15cm]
2\pp{n+\alpha+\beta+1}& n\text{ odd}
\end{array}\right.,\\[0.3cm]
    \begin{aligned}
L &=\pp{\frac{\pp{\alpha+\beta+1}x^2+\pp{c\alpha-\beta}x+c}{x^2}}\left[ R-I\right] \\
&+\left(\frac{2\pp{1-x}\pp{c+x}}{x}\right)\partial_x R .
\end{aligned}
\end{gather}\\[0.5cm]
\textbf{\textit{\large Limit Relations}}\\[0.5cm]
\textbf{\textit{\large Big q-Jacobi $\to$ Big -1 Jacobi}}\\[0.3cm]
The big -1 Jacobi polynomials are obtained from the monic big q-Jacobi polynomials by taking the parametrization $q=-e^\varepsilon$, $a=-e^{\varepsilon \alpha}$ and $b=-e^{\varepsilon \beta}$:
\begin{gather}
    \lim\limits_{\varepsilon\to0}\mathbf{P}_n\left(x;a,b,c|q\right)=\mathbf{J}_n\left(x;\alpha,\beta,c\right).
\end{gather}\\
\textbf{\textit{\large Big -1 Jacobi $\to$ Chihara}}\\[0.3cm]
The Chihara polynomials are obtained by the Christoffel transformation of the big -1 Jacobi polynomials and vice-versa via the Geronimus transformation:
\begin{gather}
    (-1)^n\pp{1-c^2}^{\frac{n}{2}}\mathbf{C}_{n}\pp{\frac{-x}{\sqrt{1-c^2}}}=\frac{\mathbf{J}_{n+1}(x)-A_{n} \mathbf{J}_{n}(x)}{x-1},\\[0.3cm] \mathbf{J}_{n}(x)=(-1)^n\pp{1-c^2}^{\frac{n}{2}}\pp{\mathbf{C}_{n}\pp{\frac{-x}{\sqrt{1-c^2}}}+\frac{C_{n}}{\sqrt{1-c^2}} \mathbf{C}_{n-1}\pp{\frac{-x}{\sqrt{1-c^2}}} },
\end{gather}
where $A_n$ and $C_n$ are given by \eqref{B-1JRRC123}. 
\begin{eqnarray}
    \frac{1}{\pp{1-c^2}^{\frac{n}{2}}}\mathbf{J}_n(x\sqrt{1-c^{2}};\alpha,\beta,c) \stackrel[\text{CT}]{\text{GT}}\longleftrightarrows \mathbf{C}_n\left(x;\frac{\beta-1}{2},\frac{\alpha+1}{2},\frac{-c}{\sqrt{1-c^{2}}}\right).
\end{eqnarray}\\
\begin{minipage}{\textwidth}\textbf{\textit{\large Continuous Bannai-Ito $\to$ Big -1 Jacobi}}\\[0.3cm]
The Big -1 Jacobi polynomial can be obtained from the Continuous Bannai-Ito polynomials using the parametrization $\beta\to\frac{\beta}{h}$, $\delta\to\frac{\delta}{h}$, the scaling $x\to2\beta x$ and taking the limit $h\to0$:
\begin{gather}
    \lim\limits_{h\to0}\left(\frac{h}{2\beta}\right)^n \mathbf{Q}_n\left(\frac{2\beta x}{h};\alpha,\frac{\beta}{h},\gamma,\frac{\delta}{h}\right) = \mathbf{J}_n\left(x;4\alpha+1,4\gamma+1,\frac{-\delta}{\beta}\right).
\end{gather}\end{minipage}\\[0.5cm]
\textbf{\textit{\large Big -1 Jacobi $\to$ Little -1 Jacobi }}\\[0.3cm]
The little -1 Jacobi polynomial can be obtained from the big -1 Jacobi polynomials by taking $c$ to 0:
\begin{gather}
    \mathbf{J}_n\pp{x;\alpha,\beta,0}=\mathbf{P}_n\pp{x;\beta,\alpha}.
\end{gather}

\subsection{Chihara}
\textbf{\textit{\large Hypergeometric Representation}}\\
\begin{flalign}
    &\mathbf{C}_{2n}\left(x;\alpha,\beta,\gamma\right)= \frac{(-1)^n\poch{\alpha+1}{n}}{\poch{n+\alpha+\beta+1}{n}}\pFq{2}{1}{-n,n+\alpha+\beta+1}{\alpha+1}{x^2-\gamma^2},& \label{CBHR1}\\
    &\mathbf{C}_{2n+1}\left(x;\alpha,\beta,\gamma\right)= \frac{(-1)^n\poch{\alpha+2}{n}}{\poch{n+\alpha+\beta+2}{n}}\pp{x-\gamma} \pFq{2}{1}{-n,n+\alpha+\beta+2}{\alpha+2}{x^2-\gamma^2}.&\label{CBHR2}
\end{flalign}
\textbf{\textit{\large Orthogonality Relation}}\\
If $\alpha>-1$ and $\beta>0$
\begin{gather}
     \int_\mathcal{C} \omega(x) \mathbf{C}_n(x)\mathbf{C}_m(x) \mathrm{d} x=h_n \delta_{n m},\\[0.3cm]
     \mathcal{C} = \bb{-\sqrt{1+\gamma^2},-|\gamma|}\bigcup\bb{|\gamma|,\sqrt{1+\gamma^2}},\\[0.3cm]
     \omega(x) = \theta(x)\pp{x+\gamma}\pp{x^2-\gamma^2}^\alpha \pp{1+\gamma^2-x^2}^{\beta},\\[0.3cm]
     \begin{array}{l}
h_{2 n}=\frac{\Gamma(n+\alpha+1) \Gamma(n+\beta+1)}{\Gamma(n+\alpha+\beta+1)} \frac{n !}{(2 n+\alpha+\beta+1)\left[(n+\alpha+\beta+1)_{n}\right]^{2}}, \\[0.3cm]
h_{2 n+1}=\frac{\Gamma(n+\alpha+2) \Gamma(n+\beta+1)}{\Gamma(n+\alpha+\beta+2)} \frac{n !}{(2 n+\alpha+\beta+2)\left[(n+\alpha+\beta+2)_{n}\right]^{2}}.
\end{array}
\end{gather}\\
\textbf{\textit{\large Normalized Recurrence Relation}}\\
\begin{gather}
    x\mathbf{C}_n(x) = \mathbf{C}_{n+1}(x) + (-1)^n\gamma\mathbf{C}_n(x) + \sigma_n \mathbf{C}_{n-1}(x),\\[0.3cm]
    \sigma_{2n}=\frac{n\pp{n+\beta}}{\pp{2n+\alpha+\beta}\pp{2n+\alpha+\beta+1}},\qquad \sigma_{2n+1}=\frac{\pp{n+\alpha+1}\pp{n+\alpha+\beta+1}}{\pp{2n+\alpha+\beta+1}\pp{2n+\alpha+\beta+2}}.
\end{gather}\\
\begin{minipage}{\textwidth}\textbf{\textit{\large Differential Equation}}
\begin{gather}
    L^{(\varepsilon)} \mathbf{C}_{n}(x)=\lambda_{n}^{(\varepsilon)} \mathbf{C}_{n}(x), \quad \left\{\begin{array}{l}
\lambda_{2n}^{(\varepsilon)}=n^2+\pp{\alpha+\beta+1}n\\[0.15cm]
\lambda_{2n+1}^{(\varepsilon)}=n^2+\pp{\alpha+\beta+2}n+\varepsilon
\end{array}\right.,\\[0.3cm]
L^{(\varepsilon)} = S(x)\partial_x^2+T(x)\partial_x R+U(x)\partial_x+V(x)\bb{I-R},\\[0.3cm]
\begin{array}{l}
S(x)=\frac{\left(x^{2}-\gamma^{2}\right)\left(x^{2}-\gamma^{2}-1\right)}{4 x^{2}}, \qquad T(x)=\frac{\gamma(x-\gamma)\left(x^{2}-\gamma^{2}-1\right)}{4 x^{3}}, \\[0.3cm]
U(x)=\frac{\gamma\left(x^{2}-\gamma^{2}-1\right)(2 \gamma-x)}{4 x^{3}}+\frac{\left(x^{2}-\gamma^{2}\right)(\alpha+\beta+3 / 2)}{2 x}-\frac{\alpha+1 / 2}{2 x}, \\[0.3cm]
V(x)=\frac{\gamma\left(x^{2}-\gamma^{2}-1\right)(x-3 \gamma / 2)}{4 x^{4}}-\frac{\left(x^{2}-\gamma^{2}\right)(\alpha+\beta+3 / 2)}{4 x^{2}}+\frac{\alpha+1 / 2}{4 x^{2}}+\epsilon \frac{x-\gamma}{2 x} .
\end{array}
\end{gather}\end{minipage}\\[0.5cm]
\textbf{\textit{\large Limit Relations}}\\[0.5cm]
\textbf{\textit{\large Big q-Jacobi $\to$ Chihara}}\\[0.3cm]
The Chihara polynomials are obtained from the monic big q-Jacobi polynomials by taking the parametrization $q=-e^\varepsilon$, $a=e^{2\varepsilon \beta}$ and $b=-e^{\varepsilon \pp{2\alpha+1}}$:
\begin{gather}
    \lim\limits_{\varepsilon\to0}\frac{1}{\pp{1-c^2}^{\frac{n}{2}}}\mathbf{P}_n\left(\sqrt{1-c^2}x;a,b,c|q\right)=\mathbf{C}_n\left(x;\alpha,\beta,\frac{-c}{\sqrt{1-c^2}}\right).
\end{gather}\\
\textbf{\textit{\large Big -1 Jacobi $\to$ Chihara}}\\[0.3cm]
The Chihara polynomials are obtained by the Christoffel transformation of the big -1 Jacobi polynomials and vice-versa via the Geronimus transformation:
\begin{gather}
    (-1)^n\pp{1-c^2}^{\frac{n}{2}}\mathbf{C}_{n}\pp{\frac{-x}{\sqrt{1-c^2}}}=\frac{\mathbf{J}_{n+1}(x)-A_{n} \mathbf{J}_{n}(x)}{x-1},\\[0.3cm] \mathbf{J}_{n}(x)=(-1)^n\pp{1-c^2}^{\frac{n}{2}}\pp{\mathbf{C}_{n}\pp{\frac{-x}{\sqrt{1-c^2}}}+\frac{C_{n}}{\sqrt{1-c^2}} \mathbf{C}_{n-1}\pp{\frac{-x}{\sqrt{1-c^2}}} },
\end{gather}
where $A_n$ and $C_n$ are given by \eqref{B-1JRRC123}. 
\begin{eqnarray}
    \frac{1}{\pp{1-c^2}^{\frac{n}{2}}}\mathbf{J}_n(x\sqrt{1-c^{2}};\alpha,\beta,c) \stackrel[\text{CT}]{\text{GT}}\longleftrightarrows \mathbf{C}_n\left(x;\frac{\beta-1}{2},\frac{\alpha+1}{2},\frac{-c}{\sqrt{1-c^{2}}}\right).
\end{eqnarray}\\
\textbf{\textit{\large Chihara $\to$ -1 Meixner-Pollaczek }}\\[0.3cm]
The -1 Meixner-Pollaczek polynomials can be obtained from the Chihara polynomials by taking $x\to\beta^{-\frac{1}{2}}x$, $\alpha\to\alpha-\frac{1}{2}$ and $\gamma\to\beta^{-\frac{1}{2}}\gamma$ and letting $\beta$ go to $\infty$:
\begin{gather}
    \lim\limits_{\beta\to\infty}\mathbf{C}_n\pp{\beta^{-\frac{1}{2}}x;\alpha-\frac{1}{2},\beta,\beta^{-\frac{1}{2}}\gamma}=\mathbf{M}_n^{(-1)}\pp{x;\alpha,\gamma}.
\end{gather}\\
\textbf{\textit{\large Chihara $\to$ Generalized Gegenbauer }}\\[0.3cm]
The generalized Gegenbauer polynomials can be obtained from the Chihara polynomials by taking $\gamma$ to 0:
\begin{gather}
    \mathbf{C}_n\pp{x;\alpha,\beta,0}=\mathbf{G}_n\pp{x;\alpha,\beta}.
\end{gather}
\subsection{Continuous -1 Hahn type 1}
\textbf{\textit{\large Hypergeometric Representation}}\\
\begin{flalign}
    &\frac{\mathbf{K}_{2n}^{\pp{1}}\left(x;\alpha,\beta,\gamma\right)}{\pp{-2i}^{2n}} = &\label{C-1H11}\end{flalign}\vspace{-0.75cm}\begin{align}\nonumber
    &\xi_{2n}\kappa_{n-1}^{(2)}\left(\frac{ix}{2}-\mathfrak{b}-\frac{1}{2}\right)\pFq{4}{3}{-n+1,n+g+2,\frac{ix}{2}+\mathfrak{b}+1,\frac{-ix}{2}+\mathfrak{b}+\frac{3}{2}}{\frac{5}{2}+\mathfrak{a}+\mathfrak{b},2+\mathfrak{b}+\mathfrak{c},2+\mathfrak{b}+\mathfrak{d}}{1} + \\
    &\kappa_n^{(1)}\pFq{4}{3}{-n,n+g+1,\frac{ix}{2}+\mathfrak{b},\frac{-ix}{2}+\mathfrak{b}+\frac{1}{2}}{\frac{3}{2}+\mathfrak{a}+\mathfrak{b},1+\mathfrak{b}+\mathfrak{c},1+\mathfrak{b}+\mathfrak{d}}{1}, \nonumber\end{align}
    \begin{flalign} &\frac{\mathbf{K}_{2n+1}^{\pp{1}}\left(x;\alpha,\beta,\gamma\right)}{\pp{-2i}^{2n+1}} =&\label{C-1H12}\end{flalign}\vspace{-0.75cm}\begin{align}\nonumber &\kappa_{n}^{(2)}\left(\frac{ix}{2}-\mathfrak{b}-\frac{1}{2}\right)\pFq{4}{3}{-n,n+g+2,\frac{ix}{2}+\mathfrak{b}+1,\frac{-ix}{2}+\mathfrak{b}+\frac{3}{2}}{\frac{5}{2}+\mathfrak{a}+\mathfrak{b},2+\mathfrak{b}+\mathfrak{c},2+\mathfrak{b}+\mathfrak{d}}{1}+\\
    &\eta_{2n+1}\kappa_n^{(1)}\pFq{4}{3}{-n,n+g+1,\frac{ix}{2}+\mathfrak{b},\frac{-ix}{2}+\mathfrak{b}+\frac{1}{2}}{\frac{3}{2}+\mathfrak{a}+\mathfrak{b},1+\mathfrak{b}+\mathfrak{c},1+\mathfrak{b}+\mathfrak{d}}{1}, \nonumber
\end{align}
where
\begin{gather}
    g=\mathfrak{a}+\mathfrak{b}+\mathfrak{c}+\mathfrak{d}+1,\qquad \mathfrak{a}=\overline{\mathfrak{d}}=\alpha+i\beta, \qquad \mathfrak{b}=\overline{\mathfrak{c}}=\gamma+i\beta,
\end{gather}
\begin{gather}
\begin{array}{ll}
\xi_{2n}=\frac{n\pp{n+\mathfrak{c}+\mathfrak{d}+\frac{1}{2}}}{\left(2n+g\right)}, & \kappa_{n}^{(1)}=\frac{\poch{\frac{3}{2}+\mathfrak{a}+\mathfrak{b}}{n} \poch{1+\mathfrak{b}+\mathfrak{c}}{n} \poch{1+\mathfrak{b}+\mathfrak{d}}{n}}{\poch{n+g+1}{n}}, \\[0.15cm]
\eta_{2n+1}=\frac{\pp{n+\mathfrak{b}+\mathfrak{c}+1}\pp{n+\mathfrak{b}+\mathfrak{d}+1}}{\pp{2n+g+1}}, &  \kappa_{n}^{(2)}=\frac{\poch{\frac{5}{2}+\mathfrak{a}+\mathfrak{b}}{n} \poch{2+\mathfrak{b}+\mathfrak{c}}{n} \poch{2+\mathfrak{b}+\mathfrak{d}}{n}}{\poch{n+g+2}{n}}.
\end{array}
\end{gather}\\
\textbf{\textit{\large Orthogonality Relation}}\\
If $\alpha,\beta,\gamma\in\mathbb{R}^+$
\begin{gather}
    \int_{-\infty}^{\infty} W(x) \mathbf{K}_{n}^{\pp{1}}(x) \mathbf{K}_{m}^{\pp{1}}(x) \mathrm{d} x=4\pi h_{0}\kappa_n \delta_{n m}, \\[0.3cm]
    W(x)=\left|\frac{\Gamma(\mathfrak{a}+i x / 2+1) \Gamma(\mathfrak{b}+i x / 2+1) \Gamma(\mathfrak{c}+i x / 2+1 / 2) \Gamma(\mathfrak{d}+i x / 2+1 / 2)}{\Gamma(1 / 2+i x)}\right|^2,\\[0.3cm]
    h_{0}=\frac{\Gamma(\mathfrak{a}+\mathfrak{b}+3 / 2) \Gamma(\mathfrak{a}+\mathfrak{c}+1) \Gamma(\mathfrak{b}+\mathfrak{c}+1) \Gamma(\mathfrak{a}+\mathfrak{d}+1) \Gamma(\mathfrak{b}+\mathfrak{d}+1) \Gamma(\mathfrak{c}+\mathfrak{d}+3 / 2)}{\Gamma(\mathfrak{a}+\mathfrak{b}+\mathfrak{c}+\mathfrak{d}+2)},\end{gather}
\begin{minipage}{\textwidth}
    \begin{flalign}
    &\kappa_{2n} = \frac{4^{2n}\Gamma\pp{n+1}\poch{2\alpha+1}{n}\poch{2\gamma+1}{n}\poch{1+\alpha+\gamma}{n}^2}{\poch{2\alpha+2\gamma+2}{2n} \poch{n+2\alpha+2\gamma+2}{n}} \bb{\prod\limits_{k=1}^n \| k+\alpha+\gamma+\frac{1}{2}+2i\beta\|}^2,&\end{flalign}
    \end{minipage}
    \begin{minipage}{\textwidth}
    \begin{flalign}
    &\kappa_{2n+1} = \frac{4^{2n+1}\Gamma\pp{n+1}\poch{2\alpha+1}{n+1}\poch{2\gamma+1}{n+1}\poch{1+\alpha+\gamma}{n+1}^2}{\poch{2\alpha+2\gamma+2}{2n+1} \poch{n+2\alpha+2\gamma+2}{n+1}}  \bb{\prod\limits_{k=1}^n \| k+\alpha+\gamma+\frac{1}{2}+2i\beta\|}^2.&
\end{flalign}\end{minipage}\\
\textbf{\textit{\large Normalized Recurrence Relation}}\\
\begin{gather}
    x\mathbf{K}_n^{(1)}(x) = \mathbf{K}_{n+1}^{(1)}(x) + b_n \mathbf{K}_n^{(1)}(x) + u_n \mathbf{K}_{n-1}^{(1)}(x),
\end{gather}
\begin{gather}
    b_{n}=\left\{\begin{array}{ll}
2 \beta-2\beta\frac{n}{(n+2 \alpha+2 \gamma+1)} & n \text { even } \\[0.3cm]
2 \beta-2\beta\frac{(n+4 \alpha+4 \gamma+3)}{n+2 \alpha+2 \gamma+2} & n \text { odd }
\end{array}\right.,\\[0.3cm]
u_{n}=\left\{\begin{array}{ll}
\frac{n(n+4 \alpha+4 \gamma+2)\|n+2\alpha+2\gamma+4i\beta+1\|^{2}}{4(n+2 \alpha+2 \gamma+1)^{2}} & n \text { even } \\[0.3cm]
\frac{(n+4 \alpha+1)(n+4 \gamma+1)\|n+2\alpha+2\gamma+1\|^{2}}{4(n+2\alpha+2 \gamma+1)^{2}} & n \text{ odd}
\end{array}\right..
\end{gather}\\
\textbf{\textit{\large Difference Equation}}\\
\begin{gather}
    L \mathbf{K}_{n}^{\pp{1}}(x)=\lambda_{n} \mathbf{K}_{n}^{\pp{1}}(x), \quad \lambda_{n}=(-1)^{n}(n+2\alpha+2\gamma+3 / 2),\\[0.3cm]
L =\mathcal{A}\left(S^{+} R-I\right)
+\overline{\mathcal{A}}\left(S^{-} R-I\right)+(2\alpha+2\gamma+3 / 2) I,\\[0.3cm]
\mathcal{A} = \frac{\pp{2 \alpha+1+i\bb{\beta- x}}\pp{2 \gamma+1+i\bb{\beta- x}}}{1-2 i x},
\end{gather}
where $\overline{\mathcal{A}}$ is the complex conjugate of $\mathcal{A}$, $S^{\pm}f(x)=f(x\pm i)$ and $Rf(x)=f(-x)$.\\[0.5cm]
\textbf{\textit{\large Limit Relations}}\\[0.3cm]
\textbf{\textit{\large Continuous Bannai-Ito $\to$ First type Continuous -1 Hahn}}\\[0.3cm]
The first type Continuous -1 Hahn polynomials can be obtained from the Continuous Bannai-Ito polynomials by a specialization:
\begin{gather}
    \mathbf{Q}_n\left(x;\alpha,\beta,\gamma,\beta\right) = \mathbf{K}_n^{(1)}\left(x;\alpha,\beta,\gamma\right).
\end{gather}\\
\begin{minipage}{\textwidth}
\textbf{\textit{\large Continuous q-Hahn $\to$ Continuous -1 Hahn type 1}}\\[0.3cm]
The Continuous -1 Hahn type 1 polynomials are obtained from the monic continuous q-Hahn polynomials by taking the parametrization $q=-e^\varepsilon$, $a=e^{\varepsilon \pp{2\alpha+1}}$, $b=e^{\varepsilon \pp{2\gamma+1}}$ and $\phi= \frac{\pi}{2}+2\varepsilon\beta$:
\begin{gather}
    \lim\limits_{\varepsilon\to0}\frac{1}{(1+q)^n}\mathbf{P}_n\left(\pp{1+q}x;a,b,\phi|q\right)=\mathbf{K}_n^{(1)}\left(x;\alpha,\beta,\gamma\right).
\end{gather}\end{minipage}\\[0.5cm]
\textbf{\textit{\large Continuous -1 Hahn $\to$ -1 Meixner-Pollaczek}}\\[0.3cm]
The -1 Meixner-Pollaczek polynomials are obtained from any type of continuous -1 Hahn polynomials by the following limit: 
\begin{gather}
    \lim\limits_{\gamma\to\infty}\frac{1}{\pp{2\gamma}^\frac{n}{2}}\mathbf{K}_{n}^{(i)}\pp{\sqrt{2\gamma}x,\frac{2\alpha-1}{4},\sqrt{\frac{\gamma}{2}}\beta,\gamma}=\mathbf{M}_n\pp{x;\alpha,\beta}.
\end{gather}\\
\textbf{\textit{\large Continuous -1 Hahn $\to$ Symmetric Bannai-Ito}}\\[0.3cm]
The symmetric Bannai-Ito polynomials are obtained from any type of continuous -1 Hahn polynomials by the specialization $\beta=0$:
\begin{gather}
    \mathbf{K}_{n}^{(i)}\pp{x,\alpha,0,\gamma}=\hat{\mathbf{S}}_n\pp{x;\alpha,\gamma}.
\end{gather}

\subsection{Continuous -1 Hahn type 2}
\textbf{\textit{\large Hypergeometric Representation}}\\
\begin{flalign}
    &\frac{\mathbf{K}_{2n}^{\pp{2}}\left(x;\alpha,\beta,\gamma\right)}{\pp{-2i}^{2n}} = &\label{C-1H21}\end{flalign}\vspace{-0.75cm}\begin{align}\nonumber
    &\xi_{2n}\kappa_{n-1}^{(2)}\left(\frac{ix}{2}-\mathfrak{b}-\frac{1}{2}\right)\pFq{4}{3}{-n+1,n+g+2,\frac{ix}{2}+\mathfrak{b}+1,\frac{-ix}{2}+\mathfrak{b}+\frac{3}{2}}{\frac{5}{2}+\mathfrak{a}+\mathfrak{b},2+\mathfrak{b}+\mathfrak{c},2+\mathfrak{b}+\mathfrak{d}}{1} + \\
    &\kappa_n^{(1)}\pFq{4}{3}{-n,n+g+1,\frac{ix}{2}+\mathfrak{b},\frac{-ix}{2}+\mathfrak{b}+\frac{1}{2}}{\frac{3}{2}+\mathfrak{a}+\mathfrak{b},1+\mathfrak{b}+\mathfrak{c},1+\mathfrak{b}+\mathfrak{d}}{1}, \nonumber\end{align}
    \begin{flalign} &\frac{\mathbf{K}_{2n+1}^{\pp{2}}\left(x;\alpha,\beta,\gamma\right)}{\pp{-2i}^{2n+1}} =&\label{C-1H22}\end{flalign}\vspace{-0.75cm}\begin{align}\nonumber &\kappa_{n}^{(2)}\left(\frac{ix}{2}-\mathfrak{b}-\frac{1}{2}\right)\pFq{4}{3}{-n,n+g+2,\frac{ix}{2}+\mathfrak{b}+1,\frac{-ix}{2}+\mathfrak{b}+\frac{3}{2}}{\frac{5}{2}+\mathfrak{a}+\mathfrak{b},2+\mathfrak{b}+\mathfrak{c},2+\mathfrak{b}+\mathfrak{d}}{1}+\\
    &\eta_{2n+1}\kappa_n^{(1)}\pFq{4}{3}{-n,n+g+1,\frac{ix}{2}+\mathfrak{b},\frac{-ix}{2}+\mathfrak{b}+\frac{1}{2}}{\frac{3}{2}+\mathfrak{a}+\mathfrak{b},1+\mathfrak{b}+\mathfrak{c},1+\mathfrak{b}+\mathfrak{d}}{1}, \nonumber
\end{align}
where
\begin{gather}
    g=2\alpha+2\gamma+1,\qquad \mathfrak{a}=\overline{\mathfrak{d}}=\alpha+i\beta, \qquad \mathfrak{b}=\overline{\mathfrak{c}}=\gamma-i\beta,
\end{gather}
\begin{gather}
\begin{array}{ll}
\xi_{2n}=\frac{n\pp{n+\mathfrak{c}+\mathfrak{d}+\frac{1}{2}}}{\left(2n+g\right)}, & \kappa_{n}^{(1)}=\frac{\poch{\frac{3}{2}+\mathfrak{a}+\mathfrak{b}}{n} \poch{1+\mathfrak{b}+\mathfrak{c}}{n} \poch{1+\mathfrak{b}+\mathfrak{d}}{n}}{\poch{n+g+1}{n}}, \\[0.15cm]
\eta_{2n+1}=\frac{\pp{n+\mathfrak{b}+\mathfrak{c}+1}\pp{n+\mathfrak{b}+\mathfrak{d}+1}}{\pp{2n+g+1}}, &  \kappa_{n}^{(2)}=\frac{\poch{\frac{5}{2}+\mathfrak{a}+\mathfrak{b}}{n} \poch{2+\mathfrak{b}+\mathfrak{c}}{n} \poch{2+\mathfrak{b}+\mathfrak{d}}{n}}{\poch{n+g+2}{n}}.
\end{array}
\end{gather}\\
\textbf{\textit{\large Orthogonality Relation}}\\
If $\alpha,\beta,\gamma\in\mathbb{R}^+$
\begin{gather}
    \int_{-\infty}^{\infty} W(x) \mathbf{K}_{n}^{\pp{2}}(x) \mathbf{K}_{m}^{\pp{2}}(x) \mathrm{d} x=4\pi h_{0}\kappa_n \delta_{n m}, \\[0.3cm]
    W(x)=\left|\frac{\Gamma(\mathfrak{a}+i x / 2+1) \Gamma(\mathfrak{b}+i x / 2+1) \Gamma(\mathfrak{c}+i x / 2+1 / 2) \Gamma(\mathfrak{d}+i x / 2+1 / 2)}{\Gamma(1 / 2+i x)}\right|^2,\\[0.3cm]
    h_{0}=\frac{\Gamma(\mathfrak{a}+\mathfrak{b}+3 / 2) \Gamma(\mathfrak{a}+\mathfrak{c}+1) \Gamma(\mathfrak{b}+\mathfrak{c}+1) \Gamma(\mathfrak{a}+\mathfrak{d}+1) \Gamma(\mathfrak{b}+\mathfrak{d}+1) \Gamma(\mathfrak{c}+\mathfrak{d}+3 / 2)}{\Gamma(\mathfrak{a}+\mathfrak{b}+\mathfrak{c}+\mathfrak{d}+2)},\end{gather}
\begin{minipage}{\textwidth}
    \begin{flalign}
    &\kappa_{2n} = \frac{4^{2n}\Gamma\pp{n+1}\poch{2\alpha+1}{n}\poch{2\gamma+1}{n}}{\poch{2\alpha+2\gamma+2}{2n} \poch{n+2\alpha+2\gamma+2}{n}} \poch{\alpha+\gamma+\frac{3}{2}}{n}^2\bb{\prod\limits_{k=1}^n \| k+\alpha+\gamma+2i\beta\|}^2,&\end{flalign}
    \end{minipage}
    \begin{minipage}{\textwidth}
    \begin{flalign}
    &\kappa_{2n+1} = \frac{4^{2n+1}\Gamma\pp{n+1}\poch{2\alpha+1}{n+1}\poch{2\gamma+1}{n+1}}{\poch{2\alpha+2\gamma+2}{2n+1} \poch{n+2\alpha+2\gamma+2}{n+1}}\poch{\alpha+\gamma+\frac{3}{2}}{n}^2\bb{\prod\limits_{k=1}^{n+1} \| k+\alpha+\gamma+2i\beta\|}^2 .&
\end{flalign}\end{minipage}\\
\textbf{\textit{\large Normalized Recurrence Relation}}\\
\begin{gather}
    x\mathbf{K}_n^{(2)}(x) = \mathbf{K}_{n+1}^{(2)}(x) + b_n \mathbf{K}_n^{(2)}(x) + u_n \mathbf{K}_{n-1}^{(2)}(x),
\end{gather}
\begin{gather}
    b_{n}=\left\{\begin{array}{ll}
2 \beta-2\beta\frac{(n+4 \alpha+2)}{(n+2 \alpha+2 \gamma+2)} & n \text { even } \\[0.3cm]
2 \beta-2\beta\frac{(n+4 \gamma+1)}{n+2 \alpha+2 \gamma+1} & n \text { odd }
\end{array}\right.,\\[0.3cm]
u_{n}=\left\{\begin{array}{ll}
\frac{n(n+4 \alpha+4 \gamma+2)\|n+2\alpha+2\gamma+1\|^{2}}{4(n+2 \alpha+2 \gamma+1)^{2}} & n \text { even } \\[0.3cm]
\frac{(n+4 \alpha+1)(n+4 \gamma+1)\|n+2\alpha+2\gamma+4i\beta+1\|^{2}}{4(n+2\alpha+2 \gamma+1)^{2}} & n \text{ odd}
\end{array}\right..
\end{gather}

\textbf{\textit{\large Difference Equation}}\\
\begin{gather}
    L \mathbf{K}_{n}^{\pp{2}}(x)=\lambda_{n} \mathbf{K}_{n}^{\pp{2}}(x), \quad \lambda_{n}=(-1)^{n}(n+2\alpha+2\gamma+3 / 2),\\[0.3cm]
L =\mathcal{A}\left(S^{+} R-I\right)
+\overline{\mathcal{A}}\left(S^{-} R-I\right)+(2\alpha+2\gamma+3 / 2) I,\\[0.3cm]
\mathcal{A} = \frac{\pp{2 \alpha+1+i\bb{\beta- x}}\pp{2 \gamma+1-i\bb{\beta+ x}}}{1-2 i x},
\end{gather}
where $\overline{\mathcal{A}}$ is the complex conjugate of $\mathcal{A}$, $S^{\pm}f(x)=f(x\pm i)$ and $Rf(x)=f(-x)$.\\[0.5cm]
\textbf{\textit{\large Limit Relations}}\\[0.5cm]
\textbf{\textit{\large Continuous Bannai-Ito $\to$ Continuous -1 Hahn Type 2}}\\[0.3cm]
The second type of Continuous -1 Hahn polynomials can be obtained from the Continuous Bannai-Ito polynomials by a specialization:
\begin{gather}
    \mathbf{Q}_n\left(x;\alpha,\beta,\gamma,-\beta\right) = \mathbf{K}_n^{(2)}\left(x;\alpha,\beta,\gamma\right).
\end{gather}\\
\textbf{\textit{\large Continuous q-Hahn $\to$ Continuous -1 Hahn type 2}}\\[0.3cm]
The Continuous -1 Hahn type 2 polynomials are obtained from the monic continuous q-Hahn polynomials by taking the parametrization $q=-e^\varepsilon$, $a=e^{\varepsilon \pp{2\alpha+1}}$, $b=-e^{\varepsilon \pp{2\gamma+1}}$ and $\phi= \frac{\pi}{2}+2\varepsilon\beta$:
\begin{gather}
    \lim\limits_{\varepsilon\to0}\frac{1}{(1+q)^n}\mathbf{P}_n\left(\pp{1+q}x;a,b,\phi|q\right)=\mathbf{K}_n^{(2)}\left(x;\alpha,\beta,\gamma\right).
\end{gather}\\
\textbf{\textit{\large Continuous -1 Hahn $\to$ -1 Meixner-Pollaczek}}\\[0.3cm]
The -1 Meixner-Pollaczek polynomials are obtained from any type of continuous -1 Hahn polynomials by the following limit: 
\begin{gather}
    \lim\limits_{\gamma\to\infty}\frac{1}{\pp{2\gamma}^\frac{n}{2}}\mathbf{K}_{n}^{(i)}\pp{\sqrt{2\gamma}x,\frac{2\alpha-1}{4},\sqrt{\frac{\gamma}{2}}\beta,\gamma}=\mathbf{M}_n\pp{x;\alpha,\beta}.
\end{gather}\\
\textbf{\textit{\large Continuous -1 Hahn $\to$ Symmetric Bannai-Ito}}\\[0.3cm]
The symmetric Bannai-Ito polynomials are obtained from any type of continuous -1 Hahn polynomials by the specialization $\beta=0$:
\begin{gather}
    \mathbf{K}_{n}^{(i)}\pp{x,\alpha,0,\gamma}=\hat{\mathbf{S}}_n\pp{x;\alpha,\gamma}.
\end{gather}
\subsection{Generalized Symmetric Bannai-Ito}
\textbf{\textit{\large Hypergeometric Representation}}\\
\begin{flalign}
    &\hat{\mathbf{I}}_{2n}\left(x;a,b,c\right)= \eta_{2n}\pFq{4}{3}{-n,n+a+b+c-1,ix,-ix}{a,b,c}{1}& \label{GSBI1},&\\[0.3cm]
    &\hat{\mathbf{I}}_{2n+1}\left(x;a,b,c\right)= \eta_{2n+1}x\pFq{4}{3}{-n,n+a+b+c,1+ix,1-ix}{1+a,1+b,1+c}{1},& \label{GSBI2}&
\end{flalign}
where
\begin{gather}
    \eta_{2n} = \frac{(-1)^n\poch{a}{n}\poch{b}{n}\poch{c}{n}}{\poch{n+a+b+c-1}{n}},\\[0.3cm] \eta_{2n+1} = \frac{(-1)^n\poch{1+a}{n}\poch{1+b}{n}\poch{1+c}{n}}{\poch{n+a+b+c}{n}}.
\end{gather}
\textbf{\textit{\large Orthogonality Relation}}\\
If Re$\pp{a,b,c}>0$, $a+b+c>1$ and non-real parameters occur in conjugate pairs, then
\begin{gather}
    \frac{1}{4\pi}\int_{-\infty}^{\infty} \omega(x) \hat{\mathbf{I}}_{n}(x)\hat{\mathbf{I}}_{m}(x) \mathrm{d} x=\kappa_n \delta_{n m}, \\[0.3cm]
    \omega(x) = \left|\frac{\Gamma\pp{ix}\Gamma\pp{a+ix}\Gamma\pp{b+ix}\Gamma\pp{c+ix}}{\Gamma\pp{2ix}}\right|^2,\\[0.3cm]
    \kappa_{n} = \frac{\Gamma\pp{n+a+b}\Gamma\pp{n+a+c}\Gamma\pp{n+b+c}\Gamma\pp{n+a}\Gamma\pp{n+b}\Gamma\pp{n+c}n!}{\Gamma\pp{2n+a+b+c}\poch{n+a+b+c-1}{n}}.
\end{gather}\\
\textbf{\textit{\large Normalized Recurrence Relation}}\\
\begin{gather}
    x\tilde{\mathbf{I}}_{n}(x)=\tilde{\mathbf{I}}_{n+1}(x)+\tau_{n} \tilde{\mathbf{I}}_{n-1}(x),\\[0.3cm]
    \begin{array}{l}
\tau_{2 n}=\frac{n\left(n+a+b-1\right)\left(n+a+c-1\right)\left(n+b+c-1\right)}{(2 n+a+b+c-2)(2 n+a+b+c-1)}, \qquad
\tau_{2 n+1}=\frac{(n+a+b+c-1)\left(n+c\right)\left(n+a\right)\left(n+b\right)}{(2 n+a+b+c-1)(2 n+a+b+c)}.
\end{array}
\end{gather}\\
\begin{minipage}{\textwidth}
\textbf{\textit{\large Difference Equation}}\\
\begin{gather}
    D_\sigma \hat{\mathbf{I}}_{n}(x) = \Lambda_{n}^{(\sigma)}\hat{\mathbf{I}}_{n}(x)\label{D0diffopGSBI},\\[0.3cm]
    \begin{array}{l}
        \Lambda_{2n}^{(\sigma)} = n^2 +\pp{a+b+c-1}n, \qquad
        \Lambda_{2n+1}^{(\sigma)} = n^2 +\pp{a+b+c}n+\sigma,
    \end{array}\\[0.3cm]
    D_{\sigma} = D_0 + \frac{\sigma}{2}\pp{I-R},\\[0.3cm]
    D_0 = B(x)S^++A(x)S^-+C(x)R-\pp{A(x)+B(x)+C(x)}I,
\end{gather}
\end{minipage}
\begin{gather}
    A(x) = \frac{\pp{ix+a}\pp{ix+b}\pp{ix+c}}{2\pp{2ix+1}},\\[0.3cm]
    B(x) = \frac{\pp{ix-a}\pp{ix-b}\pp{ix-c}}{2\pp{2ix-1}},\\[0.3cm]
    C(x) = \frac{1}{2}\pp{ab+ac+bc-x^2}-A(x)-B(x).
\end{gather}\\

\textbf{\textit{\large Limit Relations}}\\[0.3cm]
\textbf{\textit{\large Generalized Symmetric Bannai-Ito $\to$ Symmetric Bannai-Ito}}\\[0.3cm]
The symmetric Bannai-Ito polynomials can be obtained from the generalized symmetric Bannai-Ito polynomials by taking the limit $c\to \infty$:
\begin{gather}
    \lim\limits_{c\to\infty}\hat{\mathbf{I}}_{n}\pp{x,a,b,c}=\hat{\mathbf{S}}_n\pp{x;a,b}.
\end{gather}\\
\textbf{\textit{\large Generalized Symmetric Bannai-Ito $\to$ Generalized Gegenbauer}}\\[0.3cm]
The generalized Gegenbauer polynomials can be obtained from the generalized symmetric Bannai-Ito polynomials by taking the limit:
\begin{gather}
    \lim\limits_{h\to\infty}\frac{1}{h^n}\hat{\mathbf{I}}_{n}\pp{hx,\frac{\beta+1}{2}+ih,\frac{\beta+1}{2}-ih,\alpha+1}=\mathbf{G}_n\pp{x;\alpha,\beta}.
\end{gather}

\subsection{Little -1 Jacobi}
\begin{minipage}{\textwidth}
\textbf{\textit{\large Hypergeometric Representation}}\\
\begin{flalign}
    &\frac{\mathbf{P}_{2n}\left(x;\alpha,\beta\right)}{\eta_{2n}}=\label{L-1JBHR1}\pFq{2}{1}{-n,\frac{2n+\alpha+\beta+2}{2}}{\frac{\alpha+1}{2}}{x^2}+ 
    \frac{2nx}{1+\alpha}\pFq{2}{1}{1-n,\frac{2n+\alpha+\beta+2}{2}}{\frac{\alpha+3}{2}}{x^2}, &\end{flalign}
    \begin{flalign}
    &\frac{\mathbf{P}_{2n+1}\left(x;\alpha,\beta\right)}{\eta_{2n+1}}=\label{L-1JBHR2} \pFq{2}{1}{-n,\frac{n+\alpha+\beta+1}{2}}{\frac{\alpha+1}{2}}{x^2}- 
    \frac{(2n+\alpha+\beta+2)x}{1+\alpha}\pFq{2}{1}{-n,\frac{2n+\alpha+\beta+4}{2}}{\frac{\alpha+3}{2}}{x^2} ,&
\end{flalign}\end{minipage}
where
\begin{gather}
    \eta_{2n} = \frac{\poch{\frac{\alpha+1}{2}}{n}}{\poch{\frac{2n+\alpha+\beta+2}{2}}{n}}, \qquad \eta_{2n+1} = \frac{\poch{\frac{\alpha+1}{2}}{n+1}}{\poch{\frac{2n+\alpha+\beta+2}{2}}{n+1}}.
\end{gather}
\textbf{\textit{\large Orthogonality Relation}}\\
If $\alpha>0$ and $\beta>0$
\begin{gather}
     \int_{-1}^{1} \omega(x) \mathbf{P}_n(x)\mathbf{P}_m(x) \mathrm{d} x=\frac{\Gamma\pp{\frac{\alpha+1}{2}}\Gamma\pp{\frac{\beta+1}{2}}}{\Gamma\pp{\frac{\alpha}{2}+\frac{\beta}{2}+1}} \kappa_n \delta_{n m}, \\[0.3cm]
     \omega(x) = |x|^\alpha \pp{1-x^2}^{\frac{\beta-1}{2}}\pp{1+x},\\[0.3cm]
\begin{array}{l}
    \kappa_{2n} = \frac{\Gamma\pp{n+1}\poch{\frac{\alpha+1}{2}}{n}\poch{\frac{\beta+1}{2}}{n}}{\poch{1+\frac{\alpha+\beta}{2}}{2n}\poch{n+1+\frac{\alpha+\beta}{2}}{n}}, \qquad
    \kappa_{2n+1} = \frac{\Gamma\pp{n+1}\poch{\frac{\alpha+1}{2}}{n+1}\poch{\frac{\beta+1}{2}}{n+1}}{\poch{1+\frac{\alpha+\beta}{2}}{2n+1}\poch{n+1+\frac{\alpha+\beta}{2}}{n+1}} .
\end{array}
\end{gather}\\
\textbf{\textit{\large Normalized Recurrence Relation}}\\
\begin{gather}
    x\mathbf{P}_n(x) = \mathbf{P}_{n+1}(x) + \left(1-A_n-C_n\right)\mathbf{P}_n(x) + A_{n-1}C_n \mathbf{P}_{n-1}(x),
\end{gather}
\begin{minipage}{.45\textwidth}
        \begin{eqnarray}
    A_{n}=\left\{\begin{array}{ll}
\frac{n+\beta+1}{2n+\alpha+\beta+2}&n\text{ even}\\[0.15cm]
\frac{n+\alpha+\beta+1}{2n+\alpha+\beta+2}& n\text{ odd}
\end{array}\right.,\nonumber
\end{eqnarray}
    \end{minipage}%
    \begin{minipage}{0.55\textwidth}
       \begin{eqnarray}
   C_{n}=\left\{\begin{array}{ll}
\frac{n}{2n+\alpha+\beta}&n\text{ even}\\[0.15cm]
\frac{n+\alpha}{2n+\alpha+\beta}& n\text{ odd}
\end{array}\right..\label{L-1JRRC}
\end{eqnarray}
    \end{minipage}\\[0.5cm]
\begin{minipage}{\textwidth}
\textbf{\textit{\large Differential Equation}}\\
\begin{gather}
    L \mathbf{P}_n(x)=\lambda_{n} \mathbf{P}_n(x), \quad \lambda_{n}=\left\{\begin{array}{ll}
-2n&n\text{ even}\\[0.15cm]
2\pp{n+\alpha+\beta+1}& n\text{ odd}
\end{array}\right.,\\[0.3cm]
    \begin{aligned}
L =\pp{\frac{\pp{\alpha+\beta+1}x^2-\alpha x}{x^2}}\left[ R-I\right]
+2\pp{1-x}\partial_x R .
\end{aligned}
\end{gather}\end{minipage}\\[0.5cm]
\textbf{\textit{\large Limit Relations}}\\[0.5cm]
\textbf{\textit{\large Little q-Jacobi $\to$ Little -1 Jacobi}}\\[0.3cm]
The little -1 Jacobi polynomials are obtained from the monic little q-Jacobi polynomials by taking the parametrization $q=-e^\varepsilon$, $a=-e^{\varepsilon \alpha}$ and $b=-e^{\varepsilon \beta}$:
\begin{gather}
    \lim\limits_{\varepsilon\to0}\mathbf{P}_n\left(x;a,b|q\right)=\mathbf{P}_n\left(x;\alpha,\beta\right).
\end{gather}\\
\textbf{\textit{\large Big -1 Jacobi $\to$ Little -1 Jacobi }}\\[0.3cm]
The little -1 Jacobi polynomial can be obtained from the big -1 Jacobi polynomials by taking $c$ to 0:
\begin{gather}
    \mathbf{J}_n\pp{x;\alpha,\beta,0}=\mathbf{P}_n\pp{x;\beta,\alpha}.
\end{gather}\\
\textbf{\textit{\large Little -1 Jacobi $\to$ Generalized Gegenbauer}}\\[0.3cm]
The generalized Gegenbauer polynomials are obtained by the Christoffel transformation of the little -1 Jacobi polynomials and vice-versa via the Geronimus transformation:
\begin{gather}
    \mathbf{G}_{n}\pp{x}=\frac{\mathbf{P}_{n+1}(x)-A_{n} \mathbf{P}_n(x)}{x-1},\\[0.3cm] \mathbf{P}_n(x)=\mathbf{G}_{n}\pp{x}-C_{n} \mathbf{G}_{n-1}\pp{x} ,
\end{gather}
where $A_n$ and $C_n$ are given by \eqref{L-1JRRC}.
\begin{eqnarray}
    \mathbf{P}_n(x;\alpha,\beta) \stackrel[\text{CT}]{\text{GT}}\longleftrightarrows \mathbf{G}_n\left(x;\frac{\alpha-1}{2},\frac{\beta+1}{2}\right).
\end{eqnarray}\\
\begin{minipage}{\textwidth}
\textbf{\textit{\large Little -1 Jacobi $\to$ Special Little -1 Jacobi }}\\[0.3cm]
The special little -1 Jacobi polynomial can be obtained from the little -1 Jacobi polynomials by taking $\alpha$ to 0:
\begin{gather}
    \mathbf{P}_n\pp{x;0,\beta}=\mathbf{P}_n\pp{x;\beta}.
\end{gather}\end{minipage}

\subsection{Generalized Gegenbauer}
\textbf{\textit{\large Hypergeometric Representation}}\\
\begin{flalign}
    &\mathbf{G}_{2n}\left(x;\alpha,\beta\right)= \frac{(-1)^n\poch{\alpha+1}{n}}{\poch{n+\alpha+\beta+1}{n}}\pFq{2}{1}{-n,n+\alpha+\beta+1}{\alpha+1}{x^2},\label{GGBHR1}&\\
    &\mathbf{G}_{2n+1}\left(x;\alpha,\beta\right)= \frac{(-1)^n\poch{\alpha+2}{n}}{\poch{n+\alpha+\beta+2}{n}}x \pFq{2}{1}{-n,n+\alpha+\beta+2}{\alpha+2}{x^2}.\label{GGBHR2}&
\end{flalign}\\
\textbf{\textit{\large Orthogonality Relation}}\\
If $\alpha>-1$ and $\beta>0$
\begin{gather}
     \int_{-1}^{1} \omega(x) \mathbf{G}_n(x)\mathbf{G}_m(x) \mathrm{d} x=h_n \delta_{n m},\\[0.3cm]
     \omega(x) = |x|^{2\alpha+1} \pp{1-x^2}^{\beta},\\[0.3cm]
     \begin{array}{l}
h_{2 n}=\frac{\Gamma(n+\alpha+1) \Gamma(n+\beta+1)}{\Gamma(n+\alpha+\beta+1)} \frac{n !}{(2 n+\alpha+\beta+1)\left[(n+\alpha+\beta+1)_{n}\right]^{2}}, \\[0.3cm]
h_{2 n+1}=\frac{\Gamma(n+\alpha+2) \Gamma(n+\beta+1)}{\Gamma(n+\alpha+\beta+2)} \frac{n !}{(2 n+\alpha+\beta+2)\left[(n+\alpha+\beta+2)_{n}\right]^{2}}.
\end{array}
\end{gather}\\
\textbf{\textit{\large Normalized Recurrence Relation}}\\
\begin{gather}
    x\mathbf{G}_n(x) = \mathbf{G}_{n+1}(x)+ \sigma_n \mathbf{G}_{n-1}(x),\\[0.3cm]
    \sigma_{2n}=\frac{n\pp{n+\beta}}{\pp{2n+\alpha+\beta}\pp{2n+\alpha+\beta+1}},\qquad \sigma_{2n+1}=\frac{\pp{n+\alpha+1}\pp{n+\alpha+\beta+1}}{\pp{2n+\alpha+\beta+1}\pp{2n+\alpha+\beta+2}} .
\end{gather}\\
\begin{minipage}{\textwidth}
\textbf{\textit{\large Differential Equation}}\\
\begin{gather}
    L^{(\varepsilon)} \mathbf{G}_{n}(x)=\lambda_{n}^{(\varepsilon)} \mathbf{G}_{n}(x), \quad \left\{\begin{array}{l}
\lambda_{2n}^{(\varepsilon)}=n^2+\pp{\alpha+\beta+1}n\\[0.15cm]
\lambda_{2n+1}^{(\varepsilon)}=n^2+\pp{\alpha+\beta+2}n+\varepsilon
\end{array}\right.,\\[0.3cm]
L^{(\varepsilon)} = S(x)\partial_x^2+U(x)\partial_x+V(x)\bb{I-R},\\[0.3cm]
\begin{array}{l}
S(x)=\frac{x^{2}-1}{4 }, \quad
U(x)=\frac{x(\alpha+\beta+3 / 2)}{2 }-\frac{\alpha+1 / 2}{2 x} ,\quad
V(x)=-\frac{(\alpha+\beta+3 / 2)}{4 }+\frac{\alpha+1 / 2}{4 x^{2}}+ \frac{\epsilon}{2} .
\end{array}
\end{gather}\end{minipage}\\[0.5cm]
\textbf{\textit{\large Limit Relations}}\\[0.5cm]
\textbf{\textit{\large Little q-Jacobi $\to$ Generalized Gegenbauer}}\\[0.3cm]
The generalized Gegenbauer polynomials are obtained from the monic little q-Jacobi polynomials by taking the parametrization $q=-e^\varepsilon$, $b=e^{2\varepsilon \beta}$ and $a=-e^{\varepsilon \pp{2\alpha+1}}$:
\begin{gather}
    \lim\limits_{\varepsilon\to0}\mathbf{P}_n\left(x;a,b|q\right)=\mathbf{G}_n\left(x;\alpha,\beta\right).
\end{gather}\\
\textbf{\textit{\large Chihara $\to$ Generalized Gegenbauer }}\\[0.3cm]
The generalized Gegenbauer polynomials can be obtained from the Chihara polynomials by taking $\gamma$ to 0:
\begin{gather}
    \mathbf{C}_n\pp{x;\alpha,\beta,0}=\mathbf{G}_n\pp{x;\alpha,\beta}.
\end{gather}\\
\textbf{\textit{\large Little -1 Jacobi $\to$ Generalized Gegenbauer}}\\[0.3cm]
The generalized Gegenbauer polynomials are obtained by the Christoffel transformation of the little -1 Jacobi polynomials and vice-versa via the Geronimus transformation:
\begin{gather}
    \mathbf{G}_{n}\pp{x}=\frac{\mathbf{P}_{n+1}(x)-A_{n} \mathbf{P}_n(x)}{x-1},\\[0.3cm] \mathbf{P}_n(x)=\mathbf{G}_{n}\pp{x}-C_{n} \mathbf{G}_{n-1}\pp{x},
\end{gather}
where $A_n$ and $C_n$ are given by \eqref{L-1JRRC}.
\begin{eqnarray}
    \mathbf{P}_n(x;\alpha,\beta) \stackrel[\text{CT}]{\text{GT}}\longleftrightarrows \mathbf{G}_n\left(x;\frac{\alpha-1}{2},\frac{\beta+1}{2}\right).
\end{eqnarray}\\
\textbf{\textit{\large Generalized Symmetric Bannai-Ito $\to$ Generalized Gegenbauer}}\\[0.3cm]
The generalized Gegenbauer polynomials can be obtained from the generalized symmetric Bannai-Ito polynomials by taking the limit:
\begin{gather}
    \lim\limits_{h\to\infty}\frac{1}{h^n}\hat{\mathbf{I}}_{n}\pp{hx,\frac{\beta+1}{2}+ih,\frac{\beta+1}{2}-ih,\alpha+1}=\mathbf{G}_n\pp{x;\alpha,\beta}.
\end{gather}\\
\textbf{\textit{\large Generalized Gegenbauer $\to$ Generalized Hermite }}\\[0.3cm]
The Generalized Hermite polynomials can be obtained from the generalized Gegenbauer polynomials by taking $x\to\beta^{-\frac{1}{2}}x$, $\alpha\to\alpha-\frac{1}{2}$ and letting $\beta$ go to $\infty$:
\begin{gather}
    \lim\limits_{\beta\to\infty}\beta^\frac{n}{2}\mathbf{G}_n\pp{\beta^{-\frac{1}{2}}x;\alpha-\frac{1}{2},\beta}=\mathbf{H}_n\pp{x;\alpha}.
\end{gather}\\
\textbf{\textit{\large Generalized Gegenbauer $\to$ Gegenbauer }}\\[0.3cm]
The Gegenbauer polynomials can be obtained from the generalized Gegenbauer polynomials by taking $\alpha$ to $\frac{-1}{2}$:
\begin{gather}
    \mathbf{G}_n\pp{x;\frac{-1}{2},\beta-\frac{1}{2}}=\mathbf{G}_n\pp{x;\beta}.
\end{gather}
\subsection{-1 Meixner Pollaczek}
\textbf{\textit{\large Hypergeometric Representation}}\\
\begin{flalign}
    &\mathbf{M}_{2n}\left(x;\alpha,\gamma\right)= (-1)^n\poch{\alpha+\frac{1}{2}}{n}\pFq{1}{1}{-n}{\alpha+\frac{1}{2}}{x^2-\gamma^2},\label{-1MPBHR1}&\\
    &\mathbf{M}_{2n+1}\left(x;\alpha,\gamma\right)= (-1)^n\poch{\alpha+\frac{3}{2}}{n}\pp{x-\gamma} \pFq{1}{1}{-n}{\alpha+\frac{3}{2}}{x^2-\gamma^2}.\label{-1MPBHR2}&
\end{flalign}\\
\begin{minipage}{\textwidth}
\textbf{\textit{\large Orthogonality Relation}}\\
\begin{gather}
     \int_\mathcal{C} \omega(x) \mathbf{M}_n(x)\mathbf{M}_m(x) \mathrm{d} x=h_n \delta_{n m},\qquad
     \mathcal{C} = \left(-\infty,-|\gamma|\right]\bigcup\left[|\gamma|,\infty\right), \\[0.3cm]
     \omega(x) = \theta(x)\pp{x+\gamma}\pp{x^2-\gamma^2}^{\alpha-\frac{1}{2}}e^{-x^2},\\[0.3cm]
     \begin{array}{l}
h_{2 n}=n! e^{-\gamma^2}\Gamma\pp{n+\alpha+
\frac{1}{2}}, \qquad
h_{2 n+1}=n! e^{-\gamma^2}\Gamma\pp{n+\alpha+
\frac{3}{2}}.
\end{array}
\end{gather}\end{minipage}\\[0.5cm]
\textbf{\textit{\large Normalized Recurrence Relation}}\\
\begin{gather}
    x\mathbf{M}_n(x) = \mathbf{M}_{n+1}(x)+ (-1)^n\gamma \mathbf{M}_{n}(x)+ u_n \mathbf{M}_{n-1}(x),\\[0.3cm]
    u_{2n}=n,\qquad u_{2n+1}=n+\alpha+\frac{1}{2}.
\end{gather}\\
\textbf{\textit{\large Differential Equation}}\\
\begin{gather}
    L^{(\varepsilon)} \mathbf{M}_{n}(x)=\lambda_{n}^{(\varepsilon)} \mathbf{M}_{n}(x), \quad \left\{\begin{array}{l}
\lambda_{2n}^{(\varepsilon)}=n\\[0.15cm]
\lambda_{2n+1}^{(\varepsilon)}=n+\varepsilon
\end{array}\right.,\\[0.3cm]
L^{(\varepsilon)} = S(x)\partial_x^2-T(x)\partial_x R+U(x)\partial_x+V(x)\bb{I-R},\\[0.3cm]
\begin{array}{c}
S\pp{x}=\frac{\gamma^{2}-x^{2}}{4 x^{2}}, \quad T\pp{x}=\frac{\gamma(x-\gamma)}{4 x^{3}}, \\[0.3cm]
U\pp{x}=\frac{x}{2}+\frac{\gamma}{4 x^{2}}-\frac{\gamma^{2}}{2 x^{3}}-\frac{\alpha+\gamma^{2}}{2 x}, \quad V\pp{x}=\frac{3 \gamma^{2}}{8 x^{4}}-\frac{\gamma}{4 x^{3}}+\frac{\alpha+\gamma^{2}}{4 x^{2}}+\epsilon \frac{x-\gamma}{2 x}-\frac{1}{4} .
\end{array}
\end{gather}\\
\textbf{\textit{\large Limit Relations}}\\[0.3cm]
\textbf{\textit{\large q-Meixner-Pollaczek $\to$ -1 Meixner-Pollaczek}}\\[0.3cm]
The -1 Meixner-Pollaczek polynomials are obtained from the monic q-Meixner-Pollaczek polynomials by taking the parametrization $x\to\sqrt{1+q}x$, $q=-e^\varepsilon$, $a=-e^{\varepsilon \pp{\alpha+\frac{1}{2}}}$ and $\phi=\frac{\pi}{2}+\sqrt{\varepsilon}\gamma$:
\begin{gather}
    \lim\limits_{\varepsilon\to0}\mathbf{P}_n\left(x;a,\phi|q\right)=\mathbf{M}_n\left(x;\alpha,\gamma\right).
\end{gather}\\
\begin{minipage}{\textwidth}
\textbf{\textit{\large Chihara $\to$ -1 Meixner-Pollaczek }}\\[0.3cm]
The -1 Meixner-Pollaczek polynomials can be obtained from the Chihara polynomials by taking $x\to\beta^{-\frac{1}{2}}x$, $\alpha\to\alpha-\frac{1}{2}$ and $\gamma\to\beta^{-\frac{1}{2}}\gamma$ and letting $\beta$ go to $\infty$:
\begin{gather}
    \lim\limits_{\beta\to\infty}\mathbf{C}_n\pp{\beta^{-\frac{1}{2}}x;\alpha-\frac{1}{2},\beta,\beta^{-\frac{1}{2}}\gamma}=\mathbf{M}_n\pp{x;\alpha,\gamma}.
\end{gather}\end{minipage}\\[0.5cm]
\textbf{\textit{\large Continuous -1 Hahn $\to$ -1 Meixner-Pollaczek}}\\[0.3cm]
The -1 Meixner-Pollaczek polynomials are obtained from any type of continuous -1 Hahn polynomials by the following limit: 
\begin{gather}
    \lim\limits_{\gamma\to\infty}\frac{1}{\pp{2\gamma}^\frac{n}{2}}\mathbf{K}_{n}^{(i)}\pp{\sqrt{2\gamma}x,\frac{2\alpha-1}{4},\sqrt{\frac{\gamma}{2}\beta},\gamma}=\mathbf{M}_n\pp{x;\alpha,\beta}\qquad i=1,2.
\end{gather}\\
\textbf{\textit{\large -1 Meixner-Pollaczek $\to$ Generalized Hermite }}\\[0.3cm]
The generalized Hermite polynomials can be obtained from the -1 Meixner-Pollaczek polynomials by taking $\gamma$ to $0$:
\begin{gather}
    \mathbf{M}_n\pp{x;\alpha,0}=\mathbf{H}_n\pp{x;\alpha}.
\end{gather}
\subsection{Symmetric Bannai-Ito}
\textbf{\textit{\large Hypergeometric Representation}}\\
\begin{flalign}
    &\hat{\mathbf{S}}_{2n}\left(x;a,b\right)= \eta_{2n}\pFq{3}{2}{-n,ix,-ix}{a,b}{1},& \label{SBI1}&\\[0.3cm]
    &\hat{\mathbf{S}}_{2n+1}\left(x;a,b\right)= \eta_{2n+1}x\pFq{3}{2}{-n,1+ix,1-ix}{1+a,1+b}{1},& \label{SBI2}&
\end{flalign}
where
\begin{gather}
    \eta_{2n} = (-1)^n\poch{a}{n}\poch{b}{n},\qquad \eta_{2n+1} = (-1)^n\poch{1+a}{n}\poch{1+b}{n}.
\end{gather}
\textbf{\textit{\large Orthogonality Relation}}\\
If Re$\pp{a,b}>0$ and non-real parameters occur in conjugate pairs, then
\begin{gather}\label{SBIOrtUN1IN}
    \frac{1}{4\pi}\int_{-\infty}^{\infty} \omega(x) \hat{\mathbf{S}}_{n}(x)\hat{\mathbf{S}}_{m}(x) \mathrm{d} x=\kappa_n \delta_{n m}, \\[0.3cm]
    \omega(x) = \left|\frac{\Gamma\pp{ix}\Gamma\pp{a+ix}\Gamma\pp{b+ix}}{\Gamma\pp{2ix}}\right|^2,\\[0.3cm]
    \kappa_{n} = \Gamma\pp{n+a+b}\Gamma\pp{n+a}\Gamma\pp{n+b}n!.
\end{gather}
\textbf{\textit{\large Normalized Recurrence Relation}}\\
\begin{gather}
    x\tilde{\mathbf{S}}_{n}(x)=\tilde{\mathbf{S}}_{n+1}(x)+\tau_{n} \tilde{\mathbf{S}}_{n-1}(x),\\[0.3cm]
    \begin{array}{l}
\tau_{2 n}=n\left(n+a+b-1\right), \qquad
\tau_{2 n+1}=\left(n+a\right)\left(n+b\right).
\end{array}
\end{gather}\\
\begin{minipage}{\textwidth}
\textbf{\textit{\large Difference Equation}}\\
\begin{gather}
    D_\sigma \hat{\mathbf{S}}_{n}(x) = \Lambda_{n}^{(\sigma)}\hat{\mathbf{S}}_{n}(x),\label{D0diffopSBI}\\[0.3cm]
    \begin{array}{l}
        \Lambda_{2n}^{(\sigma)} = n ,\qquad
        \Lambda_{2n+1}^{(\sigma)} = n+\sigma,
    \end{array}\\[0.3cm]
    D_{\sigma} = D_0 + \frac{\sigma}{2}\pp{I-R},\\[0.3cm]
    D_0 = B(x)S^++A(x)S^-+C(x)R-\pp{A(x)+B(x)+C(x)}I,
\end{gather}
\begin{gather}
    A(x) = \frac{\pp{ix+a}\pp{ix+b}}{2\pp{1+2ix}},\\[0.3cm]
    B(x) = \frac{\pp{ix-a}\pp{ix-b}}{2\pp{1-2ix}},\\[0.3cm]
    C(x) = \frac{a+b}{2}-A(x)-B(x).
\end{gather}\end{minipage}\\[0.5cm]
\begin{minipage}{\textwidth}
\textbf{\textit{\large Limit Relations}}\\[0.5cm]
\textbf{\textit{\large Generalized Symmetric Bannai-Ito $\to$ Symmetric Bannai-Ito}}\\[0.3cm]
The symmetric Bannai-Ito polynomials can be obtained from the generalized symmetric Bannai-Ito polynomials by taking the limit:
\begin{gather}
    \lim\limits_{c\to\infty}\hat{\mathbf{I}}_{n}\pp{x,a,b,c}=\hat{\mathbf{S}}_n\pp{x;a,b}.
\end{gather}\end{minipage}\\[0.5cm]
\textbf{\textit{\large Continuous -1 Hahn $\to$ Symmetric Bannai-Ito}}\\[0.3cm]
The symmetric Bannai-Ito polynomials are obtained from any type of continuous -1 Hahn polynomials by the specialization $\beta=0$:
\begin{gather}
    \mathbf{K}_{n}^{(i)}\pp{x,\alpha,0,\gamma}=\hat{\mathbf{S}}_n\pp{x;\alpha,\gamma}.
\end{gather}\\
\textbf{\textit{\large Symmetric Bannai-Ito $\to$ Generalized Hermite }}\\[0.3cm]
The generalized Hermite polynomials can be obtained from the Symmetric Bannai-Ito polynomials by taking $x\to\sqrt{\beta}x$, renormalizing and letting $b$ go to $\infty$:
\begin{gather}
    \lim\limits_{b\to\infty}\frac{1}{b^{\frac{n}{2}}}\hat{\mathbf{S}}_n\pp{\sqrt{b}x;\alpha+\frac{1}{2},b}=\mathbf{H}_n\pp{x;\alpha}.
\end{gather}

\subsection{Special Little -1 Jacobi}
\textbf{\textit{\large Hypergeometric Representation}}\\
\begin{flalign}
    &\frac{\mathbf{P}_{2n}\left(x;\alpha\right)}{\eta_{2n}}= \pFq{2}{1}{-n,\frac{2n+\alpha+2}{2}}{\frac{1}{2}}{x^2}+ \label{SL-1JBHR1}
    2nx\pFq{2}{1}{1-n,\frac{2n+\alpha+2}{2}}{\frac{3}{2}}{x^2}  ,&\\
    &\frac{\mathbf{P}_{2n+1}\left(x;\alpha\right)}{\eta_{2n+1}}=\pFq{2}{1}{-n,\frac{2n+\alpha+2}{2}}{\frac{1}{2}}{x^2}- \label{SL-1JBHR2}
    (2n+\alpha+2)x\pFq{2}{1}{-n,\frac{n+\alpha+3}{2}}{\frac{3}{2}}{x^2},&
\end{flalign}
where
\begin{gather}
    k_{2n} = \frac{\poch{\frac{1}{2}}{n}}{\poch{\frac{2n+\alpha+2}{2}}{n}}, \qquad k_{2n+1}=\frac{\poch{\frac{1}{2}}{n+1}}{\poch{\frac{2n+\alpha+2}{2}}{n+1}}.
\end{gather}\\
\begin{minipage}{\textwidth}
\textbf{\textit{\large Orthogonality Relation}}\\
If $\alpha>0$
\begin{gather}
     \int_{-1}^{1} \omega(x) \mathbf{P}_n(x)\mathbf{P}_m(x) \mathrm{d} x= \pp{\frac{\sqrt{\pi}\Gamma\pp{\frac{\alpha+1}{2}}}{2^{2n}}} \pp{\frac{\Gamma\pp{n+1}\Gamma\pp{n+1+\alpha}}{\Gamma\pp{n+1+\frac{\alpha}{2}}^2}}\frac{\Gamma\pp{1+\frac{\alpha}{2}}}{\Gamma\pp{1+\alpha}} \delta_{n m}, \\[0.3cm]
     \omega(x) =  \pp{1-x^2}^{\frac{\alpha-1}{2}}\pp{1+x}.
\end{gather}\end{minipage}\\[0.5cm]
\textbf{\textit{\large Normalized Recurrence Relation}}\\
\begin{gather}
    x\mathbf{P}_n(x) = \mathbf{P}_{n+1}(x) + \left(1-A_n-C_n\right)\mathbf{P}_n(x) + A_{n-1}C_n \mathbf{P}_{n-1}(x),\\[0.3cm]
    A_{n}=\frac{n+\alpha+1}{2n+\alpha+2},\qquad
   C_{n}=\frac{n}{2n+\alpha}.\label{SL-1JRRC}
\end{gather}\\
\textbf{\textit{\large Differential Equation}}\\
\begin{gather}
    L \mathbf{P}_n(x)=\lambda_{n} \mathbf{P}_n(x), \quad \lambda_{n}=\left\{\begin{array}{ll}
-2n&n\text{ even}\\[0.15cm]
2\pp{n+\alpha+1}& n\text{ odd}
\end{array}\right.,\\[0.3cm]
    \begin{aligned}
L =\pp{\alpha+1}\left[ R-I\right]
+2\pp{1-x}\partial_x R .
\end{aligned}
\end{gather}\\
\begin{minipage}{\textwidth}
\textbf{\textit{\large Limit Relations}}\\[0.5cm]
\textbf{\textit{\large Little -1 Jacobi $\to$ Special Little -1 Jacobi }}\\[0.3cm]
The special little -1 Jacobi polynomial can be obtained from the little -1 Jacobi polynomials by taking $\alpha$ to 0:
\begin{gather}
    \mathbf{P}_n\pp{x;0,\beta}=\mathbf{P}_n\pp{x;\beta}.
\end{gather}\end{minipage}\\[0.5cm]
\textbf{\textit{\large Special Little -1 Jacobi $\to$ Gegenbauer}}\\[0.3cm]
The Gegenbauer polynomials are obtained by the Christoffel transformation of the special little -1 Jacobi polynomials and vice-versa via the Geronimus transformation:
\begin{gather}
    \mathbf{G}_{n}\pp{x}=\frac{\mathbf{P}_{n+1}(x)-A_{n} \mathbf{P}_n(x)}{x-1},\\[0.3cm] \mathbf{P}_n(x)=\mathbf{G}_{n}\pp{x}-C_{n} \mathbf{G}_{n-1}\pp{x}, 
\end{gather}
where $A_n$ and $C_n$ are given by \eqref{SL-1JRRC}.
\begin{eqnarray}
    \mathbf{P}_n(x;\alpha) \stackrel[\text{CT}]{\text{GT}}\longleftrightarrows \mathbf{G}_n\left(x;\frac{\alpha+2}{2}\right).
\end{eqnarray}

\subsection{Gegenbauer}
\textbf{\textit{\large Hypergeometric Representation}}\\
\begin{flalign}
    &\mathbf{G}_{2n}\left(x;\alpha\right)= \frac{(-1)^n\poch{\frac{1}{2}}{n}}{\poch{n+\alpha}{n}}\pFq{2}{1}{-n,n+\alpha}{\frac{1}{2}}{x^2},\label{GBHR1}&\\
    &\mathbf{G}_{2n+1}\left(x;\alpha\right)= \frac{(-1)^n\poch{\frac{3}{2}}{n}}{\poch{n+\alpha+1}{n}}x \pFq{2}{1}{-n,n+\alpha+1}{\frac{3}{2}}{x^2}.\label{GBHR2}&
\end{flalign}\\
\begin{minipage}{\textwidth}
\textbf{\textit{\large Orthogonality Relation}}\\
\begin{gather}
     \int_\mathcal{C} \omega(x) \mathbf{G}_n(x)\mathbf{G}_m(x) \mathrm{d} x=h_n \delta_{n m},\qquad
     \mathcal{C} = \bb{-1,1},\\[0.3cm]
     \omega(x) = \pp{1-x^2}^{\alpha-\frac{1}{2}},\end{gather}\end{minipage}\\[0.3cm]
     \begin{gather}
     \begin{array}{l}
h_{2 n}=\frac{\Gamma(n+\frac{1}{2}) \Gamma(n+\alpha+\frac{1}{2})}{\Gamma(n+\alpha)} \frac{n !}{(2 n+\alpha)\left[(n+\alpha)_{n}\right]^{2}}, \\[0.3cm]
h_{2 n+1}=\frac{\Gamma(n+\frac{3}{2}) \Gamma(n+\alpha+\frac{1}{2})}{\Gamma(n+\alpha+1)} \frac{n !}{(2 n+\alpha+1)\left[(n+\alpha+1)_{n}\right]^{2}}.
\end{array}
\end{gather}\\
\textbf{\textit{\large Normalized Recurrence Relation}}\\
\begin{gather}
    x\mathbf{G}_n(x) = \mathbf{G}_{n+1}(x)+ \sigma_n \mathbf{G}_{n-1}(x),\\[0.3cm]
    \sigma_{n}=\frac{n\pp{n+2\alpha-1}}{\pp{2n+2\alpha-2}\pp{2n+2\alpha}}.
\end{gather}\\
\textbf{\textit{\large Differential Equation}}\\
\begin{gather}
    L^{(\varepsilon)} \mathbf{G}_{n}(x)=\lambda_{n}^{(\varepsilon)} \mathbf{G}_{n}(x), \quad \left\{\begin{array}{l}
\lambda_{2n}^{(\varepsilon)}=n^2+\alpha n\\[0.15cm]
\lambda_{2n+1}^{(\varepsilon)}=n^2+\pp{\alpha+1}n+\varepsilon
\end{array}\right.,\\[0.3cm]
L^{(\varepsilon)} = S(x)\partial_x^2+U(x)\partial_x+V(x)\bb{I-R},\\[0.3cm]
\begin{array}{l}
S(x)=\frac{x^{2}-1}{4 } ,\quad
U(x)=\frac{x(\alpha+1 / 2)}{2 }, \quad
V(x)=-\frac{(\alpha+1 / 2)}{4 }+ \frac{\epsilon}{2} .
\end{array}
\end{gather}\\
\textbf{\textit{\large Limit Relations}}\\[0.5cm]
\textbf{\textit{\large Generalized Gegenbauer $\to$ Gegenbauer }}\\[0.3cm]
The Gegenbauer polynomials can be obtained from the generalized Gegenbauer polynomials by taking $\alpha$ to $\frac{-1}{2}$:
\begin{gather}
    \mathbf{G}_n\pp{x;\frac{-1}{2},\beta-\frac{1}{2}}=\mathbf{G}_n\pp{x;\beta}.
\end{gather}\\
\textbf{\textit{\large Special Little -1 Jacobi $\to$ Gegenbauer}}\\[0.3cm]
The Gegenbauer polynomials are obtained by the Christoffel transformation of the special little -1 Jacobi polynomials and vice-versa via the Geronimus transformation:
\begin{gather}
    \mathbf{G}_{n}\pp{x}=\frac{\mathbf{P}_{n+1}(x)-A_{n} \mathbf{P}_n(x)}{x-1},\\[0.3cm] \mathbf{P}_n(x)=\mathbf{G}_{n}\pp{x}-C_{n} \mathbf{G}_{n-1}\pp{x},
\end{gather}
where $A_n$ and $C_n$ are given by \eqref{SL-1JRRC}.
\begin{eqnarray}
    \mathbf{P}_n(x;\alpha) \stackrel[\text{CT}]{\text{GT}}\longleftrightarrows \mathbf{G}_n\left(x;\frac{\alpha+2}{2}\right).
\end{eqnarray}\\
\begin{minipage}{\textwidth}
\textbf{\textit{\large Gegenbauer $\to$ Hermite }}\\[0.3cm]
The Hermite polynomials can be obtained from the Gegenbauer polynomials by taking $x\to\beta^{-\frac{1}{2}}x$ and letting $\beta$ go to $\infty$:
\begin{gather}
    \lim\limits_{\alpha\to\infty}\alpha^\frac{n}{2}\mathbf{G}_n\pp{\alpha^{-\frac{1}{2}}x;\alpha}=\mathbf{H}_n\pp{x}.
\end{gather}\end{minipage}
\subsection{Generalized Hermite}
\textbf{\textit{\large Hypergeometric Representation}}\\
\begin{flalign}
    &\mathbf{H}_{2n}\left(x;\alpha\right)= (-1)^n\poch{\alpha+\frac{1}{2}}{n}\pFq{1}{1}{-n}{\alpha+\frac{1}{2}}{x^2},\label{GHBHR1}&\\
    &\mathbf{H}_{2n+1}\left(x;\alpha\right)= (-1)^n\poch{\alpha+\frac{3}{2}}{n}x \pFq{1}{1}{-n}{\alpha+\frac{3}{2}}{x^2}.\label{GHBHR2}
\end{flalign}\\
\textbf{\textit{\large Orthogonality Relation}}\\
\begin{gather}
     \int_\mathcal{C} \omega(x) \mathbf{H}_{n}(x)\mathbf{H}_{m}(x) \mathrm{d} x=h_n \delta_{n m},\qquad
     \mathcal{C} = \left(-\infty,\infty\right) ,\\[0.3cm]
     \omega(x) = |x|^{2\alpha}e^{-x^2},\\[0.3cm]
     \begin{array}{l}
h_{2 n}=n! \Gamma\pp{n+\alpha+
\frac{1}{2}}, \qquad
h_{2 n+1}=n! \Gamma\pp{n+\alpha+
\frac{3}{2}}.
\end{array}
\end{gather}\\
\textbf{\textit{\large Normalized Recurrence Relation}}\\
\begin{gather}
    x\mathbf{H}_{n}(x) = \mathbf{H}_{n+1}(x)+ u_n \mathbf{H}_{n-1}(x),\\[0.3cm]
    u_{2n}=n,\qquad u_{2n+1}=n+\alpha+\frac{1}{2}.
\end{gather}\\
\begin{minipage}{\textwidth}
\textbf{\textit{\large Differential Equation}}\\
\begin{gather}
    L^{(\varepsilon)} \mathbf{H}_{n}(x)=\lambda_{n}^{(\varepsilon)} \mathbf{H}_{n}(x), \quad \left\{\begin{array}{l}
\lambda_{2n}^{(\varepsilon)}=n\\[0.15cm]
\lambda_{2n+1}^{(\varepsilon)}=n+\varepsilon
\end{array}\right.,\\[0.3cm]
L^{(\varepsilon)} = S(x)\partial_x^2+U(x)\partial_x+V(x)\bb{I-R},\\[0.3cm]
\begin{array}{l}
S\pp{x}=\frac{-1}{4},\qquad U\pp{x}=\frac{x}{2}-\frac{\alpha}{2 x},\qquad
V\pp{x}=\frac{\alpha}{4 x^{2}}+\frac{\epsilon}{2} -\frac{1}{4} .
\end{array}
\end{gather}\end{minipage}\\[0.5cm]
\textbf{\textit{\large Limit Relations}}\\[0.5cm]
\textbf{\textit{\large Generalized Gegenbauer $\to$ Generalized Hermite }}\\[0.3cm]
The Generalized Hermite polynomials can be obtained from the generalized Gegenbauer polynomials by taking $x\to\beta^{-\frac{1}{2}}x$, $\alpha\to\alpha-\frac{1}{2}$ and letting $\beta$ go to $\infty$:
\begin{gather}
    \lim\limits_{\beta\to\infty}\beta^\frac{n}{2}\mathbf{G}_n\pp{\beta^{-\frac{1}{2}}x;\alpha-\frac{1}{2},\beta}=\mathbf{H}_n\pp{x;\alpha}.
\end{gather}\\
\textbf{\textit{\large -1 Meixner-Pollaczek $\to$ Generalized Hermite }}\\[0.3cm]
The generalized Hermite polynomials can be obtained from the -1 Meixner-Pollaczek polynomials by taking $\gamma$ to $0$:
\begin{gather}
    \mathbf{M}_n\pp{x;\alpha,0}=\mathbf{H}_n\pp{x;\alpha}.
\end{gather}\\
\textbf{\textit{\large Symmetric Bannai-Ito $\to$ Generalized Hermite }}\\[0.3cm]
The generalized Hermite polynomials can be obtained from the Symmetric Bannai-Ito polynomials by taking $x\to\sqrt{\beta}x$, renormalizing and letting $b$ go to $\infty$:
\begin{gather}
    \lim\limits_{b\to\infty}\frac{1}{b^{\frac{n}{2}}}\hat{\mathbf{S}}_n\pp{\sqrt{b}x;\alpha+\frac{1}{2},b}=\mathbf{H}_n\pp{x;\alpha}.
\end{gather}\\
\begin{minipage}{\textwidth}
\textbf{\textit{\large Generalized Hermite $\to$ Hermite }}\\[0.3cm]
The Hermite polynomials can be obtained from the generalized Hermite  polynomials by taking $\alpha$ to $0$:
\begin{gather}
    \mathbf{H}_n\pp{x;0}=\mathbf{H}_n\pp{x}.
\end{gather}\end{minipage}
\subsection{Hermite}
\textbf{\textit{\large Hypergeometric Representation}}\\
\begin{flalign}
    &\mathbf{H}_{2n}\left(x\right)= (-1)^n\poch{\frac{1}{2}}{n}\pFq{1}{1}{-n}{\frac{1}{2}}{x^2},\label{HBHR1}&\\
    &\mathbf{H}_{2n+1}\left(x\right)= (-1)^n\poch{\frac{3}{2}}{n}x \pFq{1}{1}{-n}{\frac{3}{2}}{x^2}.\label{HBHR2}&
\end{flalign}\\
\textbf{\textit{\large Orthogonality Relation}}\\
\begin{gather}
     \int_\mathcal{C} \omega(x) \mathbf{H}_{n}(x)\mathbf{H}_{m}(x) \mathrm{d} x=h_n \delta_{n m},\qquad
     \mathcal{C} = \left(-\infty,\infty\right), \\[0.3cm]
     \omega(x) = e^{-x^2},\\[0.3cm]
     \begin{array}{l}
h_{2 n}=n! \Gamma\pp{n+
\frac{1}{2}}, \qquad
h_{2 n+1}=n! \Gamma\pp{n+
\frac{3}{2}}.
\end{array}
\end{gather}\\
\textbf{\textit{\large Normalized Recurrence Relation}}\\
\begin{gather}
    x\mathbf{H}_{n}(x) = \mathbf{H}_{n+1}(x)+ u_n \mathbf{H}_{n-1}(x),\\[0.3cm]
    u_{n}=\frac{n}{2}.
\end{gather}\\
\begin{minipage}{\textwidth}
\textbf{\textit{\large Differential Equation}}\\
\begin{gather}
    L^{(\varepsilon)} \mathbf{H}_{n}(x)=\lambda_{n}^{(\varepsilon)} \mathbf{H}_{n}(x), \quad \left\{\begin{array}{l}
\lambda_{2n}^{(\varepsilon)}=n\\[0.15cm]
\lambda_{2n+1}^{(\varepsilon)}=n+\varepsilon
\end{array}\right.,\\[0.3cm]
L^{(\varepsilon)} = S(x)\partial_x^2+U(x)\partial_x+V(x)\bb{I-R},\\[0.3cm]
\begin{array}{l}
S\pp{x}=\frac{-1}{4},\qquad U\pp{x}=\frac{x}{2},\qquad
V\pp{x}=\frac{\epsilon}{2} -\frac{1}{4} .
\end{array}
\end{gather}\end{minipage}\\[0.5cm]
\begin{minipage}{\textwidth}
\textbf{\textit{\large Limit Relations}}\\[0.5cm]
\textbf{\textit{\large Gegenbauer $\to$ Hermite }}\\[0.3cm]
The Hermite polynomials can be obtained from the Gegenbauer polynomials by taking $x\to\beta^{-\frac{1}{2}}x$ and letting $\beta$ go to $\infty$:
\begin{gather}
    \lim\limits_{\beta\to\infty}\beta^\frac{n}{2}\mathbf{G}_n\pp{\beta^{-\frac{1}{2}}x;\beta}=\mathbf{H}_n\pp{x}.
\end{gather}\end{minipage}\\[0.5cm]
\textbf{\textit{\large Generalized Hermite $\to$ Hermite }}\\[0.3cm]
The Hermite polynomials can be obtained from the generalized Hermite  polynomials by taking $\alpha$ to $0$:
\begin{gather}
    \mathbf{H}_n\pp{x;0}=\mathbf{H}_n\pp{x}.
\end{gather}
\bibliographystyle{myabbrv}
\bibliography{biblioC_1HOP}

\end{document}